\documentclass[draftcls,onecolumn,12pt]{IEEEtran}
\usepackage{amssymb,times}
\usepackage{bbm}
\usepackage{adjustbox}
\usepackage{lipsum}%
\usepackage{cite, epsfig}
\usepackage{verbatim, amsfonts,enumerate}     
\usepackage{bm,array}
\usepackage{graphicx}
\usepackage[cmex10]{amsmath}
\usepackage{epstopdf}
\usepackage{subcaption}
\usepackage{amsthm}
\usepackage{algorithm,algorithmic}
\usepackage{color,url}
\usepackage[table]{xcolor}
\usepackage[normalem]{ulem}

\theoremstyle{plain}
\newtheorem{thm}{Theorem}
\newtheorem{lem}{Lemma}
\def\D{\mathcal{D}}
\def\U{\mathbf{U}}
\def\V{\mathbf{V}}
\def\Sig{\boldsymbol{\Sigma}}

\def\v{\mathbf{v}}
\def\u{\mathbf{u}}
\def\x{\mathbf{x}}
\def\e{\mathbf{e}}
\def\y{\mathbf{y}}
\def\R{\Rn}
\def\t{\mathbf{t}}

\def\cX {\mathcal{X}}
\def\B{\mathbf{B}}
\def\R{\mathbf{R}}
\def\L{\mathbf{L}}
\def\M{\mathbf{M}}
\def\Rn{\mathbb{R}}

\def \Ob {\boldsymbol{\Omega}}
\def\Z{\mathbf{Z}}
\def\N{\mathbb{N}}
\def\A{\mathbf{A}}
\def\S{\mathbf{S}}

\def\z{\mathbf{z}}
\def\w{\mathbf{w}}
\def\m{\boldsymbol{m}}
\def\p{\boldsymbol{p}}

\def\xks{\x_k^\star}
\def\H{\mathbf{H}}
\def \nt {\tilde{\nabla}}
\def\s{\mathbf{s}}
\def\cS{\mathcal{S}}
\def\cL{\mathcal{L}}
\def\I{\mathcal{I}}
\def\Y{\mathbf{Y}}

\def\psib{\boldsymbol{\psi}}
\def\bt{\boldsymbol{\beta}}
\def\EE{\mathbb{E}}
\def\prox{\textbf{prox}}

\providecommand{\abs}[1]{\left|#1\right|}
\providecommand{\norm}[1]{\left\|#1\right\|}
\providecommand{\ip}[1]{\langle#1\rangle}

\providecommand{\vect}[1]{\text{vec}\left(#1\right)}

\providecommand{\Ex}[1]{\mathbb{E}\left[#1\right]}
\providecommand{\pk}[1]{\prox_{g_k}^\alpha\!\!\left(#1\right)}

\newcommand{\ind}{1\hspace{-1.6mm}1}	

\theoremstyle{plain}
\theoremstyle{remark}

\graphicspath{{img/}}
	\title{Online Learning with Inexact Proximal Online Gradient Descent Algorithms
}
\author{Rishabh Dixit,	Amrit Singh Bedi,~\IEEEmembership{Student Member,~IEEE}, Ruchi Tripathi,~\IEEEmembership{Student Member,~IEEE},  
	and~Ketan~Rajawat,~\IEEEmembership{Member,~IEEE}
	\thanks{
		The authors are with the Department of Electrical Engineering,
		Indian Institute of Technology Kanpur, Kanpur 208016, India (e-mail:
		rishabd93@gmail.com; amritbd@iitk.ac.in; ruchi@iitk.ac.in; ketan@iitk.ac.in). A part \cite{asilomar} of this work is submitted to 52nd Asilomar Conf. on Signals, Systems, and Computers, Pacific Grove, CA, Nov. 2018. }\vspace{-0mm}}
	\allowdisplaybreaks
\begin{document}	

	\maketitle

\begin{abstract}
We consider non-differentiable dynamic optimization problems such as those arising in robotics and subspace tracking. Given the computational constraints and the time-varying nature of the problem, a low-complexity algorithm is desirable, while the accuracy of the solution may only increase slowly over time. We put forth the proximal online gradient descent (OGD) algorithm for tracking the optimum of a composite objective function comprising of a differentiable loss function and a non-differentiable regularizer. An online learning framework is considered and the gradient of the loss function is allowed to be erroneous. Both, the gradient error as well as the dynamics of the function optimum or target are adversarial and the performance of the inexact proximal OGD is characterized in terms of its dynamic regret, expressed in terms of the cumulative error and path length of the target. The proposed inexact proximal OGD is generalized for application to large-scale problems where the loss function has a finite sum structure. In such cases, evaluation of the full gradient may not be viable and a variance reduced version is proposed that allows the component functions to be sub-sampled. The efficacy of the proposed algorithms is tested on the problem of formation control in robotics and on the dynamic foreground-background separation problem in video. 
\end{abstract}


	\section{ Introduction}
	Time-varying optimization problems arise in a number of signal processing disciplines, such as target tracking, robotics, dynamic learning, and cyber-physical systems \cite{derenick2009optimal,mokhtari2016online,6293884}. Different from the classical online learning settings that focus on processing data in a sequential or incremental fashion in order to solve a static problem, dynamic optimization algorithms seek to track the target variables over time. The robust subspace tracking problem for instance, arising in video foreground-background separation, speech enhancement, network monitoring, and dynamic MRI, entails identifying a low-rank and a sparse component in the measurements, that arrive sequentially over time. Such streaming data is often handled via stochastic incremental or online algorithms that are capable of learning the underlying subspace by processing only a few measurements at every iteration \cite{feng2013online, lois2015online, bouwmans2014robust, he2011online}. These algorithms however approach the problem from a static perspective, producing a sequence of increasingly accurate estimates of the actual subspace \cite{reprocs, wang2012probabilistic}. In contrast, the underlying optimization variables in most such problems may generally be time-varying, e.g., due to a slowly changing background in a video. Existing algorithms utilize heuristics to track these time-varying parameters, e.g., using moving windows, and are often not amenable to performance guarantees. 
	
	Online convex optimization provides the tools to handle such dynamic problems in real-time settings. Online convex optimization algorithms are explicitly designed to have low run-time complexity while incurring a bearable loss or regret \cite{zinkevich2003online}. While initially motivated for static learning problems, online algorithms for dynamic scenarios have recently been developed and successfully implemented for a wide range of applications \cite{hall2015online, zhang2016improved, simonetto2015non, shahrampour2018distributed,pmlr-v38-jadbabaie15}. Different from the classical notion of 'static' regret, that measures the cumulative error against a static adversary, the performance of such algorithms is measured against that of a time-varying adversary that models the dynamically changing optimum. The resulting \emph{dynamic} regret always upper bounds the static regret and is characterized in terms of the cumulative variations in the problem parameters.
	
	The idea of dynamic regret is still quite nascent and has only been considered for a small class of convex and differentiable objective functions. However, a number of signal processing problems, such as robust subspace tracking, have non-differentiable regularizers to encourage low-rank and/or sparse solutions. To this end, this work introduces the proximal online gradient descent (OGD) class of algorithms, whose static variants have been widely studied. A generic analysis is provided where the gradient calculations are allowed to be erroneous and the dynamic regret is characterized in terms of the cumulative problem variations and the cumulative mean squared error. The results provided here also improve upon the bounds for non-differentiable problems in  \cite{hall2015online} while also generalizing earlier results for self-concordant functions \cite{zhang2016improved}.
	
	Further challenges arise in large-scale problems where evaluating the full gradient at every iteration may be excessively costly. The cost function in such problems is often expressible as a sum of several component functions, each depending only on a subset of measurements. In order to solve such problems in real-time it becomes necessary to trade-off solution accuracy with speed, and sample only a single or few component functions at every time. However, direct application of such a sampling process would incur a large error in the gradient that may cause the iterates to diverge. Towards this end, we propose variance reduced algorithms that utilize the idea of averaging gradient from \cite{defazio2014saga} and work in time varying scenarios. The proposed variance reduced online algorithm is seminal and provides an entirely new subsampling approach to handling big streaming data. Numerical tests on the dynamic video foreground-background separation problem showcase that the proposed algorithm is computationally efficient and capable of separating dynamic background and foreground. Before proceeding with the detailed description of the proposed algorithms, we briefly review the related work.


	
	\subsection{Related work and contributions}
Online convex optimization and the notion of regret was first introduced in \cite{zinkevich2003online} and has been widely applied to machine learning \cite{zinkevich2003online,xiao2014proximal,Duchi2011}. As an example, the adaptive subgradient methods developed in \cite{Duchi2011} improve upon the regret bounds first proposed in \cite{zinkevich2003online}. However in  more general settings, such as when the optimal itself varies with time, static regret fails as a robust measure of cumulative error. For such cases, it becomes necessary to characterize the performance of the online algorithm via dynamic regret, evaluated against a time varying comparator \cite{besbes2015non, bedi2018tracking}. 

This problem of tracking time-varying adversaries has been well-studied and the dynamic regret bounds for a number of settings are well-known \cite{hall2015online, zhang2016improved, simonetto2015non, shahrampour2018distributed,pmlr-v38-jadbabaie15, bedi2018tracking}. Most of these works however focus on the differentiable case and consequently, cannot handle standard non-differentiable regularizers such as  $\ell_1$-norm or nuclear norm penalties. Dynamic regret bounds for general non-differentiable functions were first developed in \cite{hall2015online}. Different from \cite{hall2015online}, we consider a differentiable loss function but a non-differentiable regularizer, allowing the development of tighter dynamic regret bounds. A running ADMM algorithm was also proposed in \cite{simonetto2015non} but no theoretical guarantees were provided. Subsequently, a class of proximal algorithms were proposed in \cite{simonetto2017time, bernstein2018asynchronous} and analyzed from the perspective of time-varying optimization. Different from these works however, we develop dynamic regret bounds that require a different set of assumptions. Additionally, the bounds developed here allow inexact gradients and readily extend to big-data settings where the finite sum structure of the objective function prevents us from evaluating the full gradient at each time step. It is remarked that non-smooth but differentiable functions are also allowed within the framework considered in \cite{zhang2016improved}. As such, the approach in \cite{zhang2016improved} is not applicable to the case when the regularizer is non-differentiable. 

The proof in the present work builds upon the analysis of proximal and incremental gradient methods applied to static problems in \cite{bertsekas2000gradient}. In large-scale settings, acceleration methods such as averaging gradient \cite{defazio2014saga} and variance reduction \cite{xiao2014proximal} have been widely used. The present work applies these methods to the dynamic setting yielding a novel approach to tracking with subsampling.

The problem of separating the foreground and background of a video have been well-studied. The state-of-the-art low-complexity online algorithms include  GRASTA\cite{he2012incremental}, ReProCS \cite{reprocs} and OPRMF \cite{wang2012probabilistic}. Of these, GRASTA is a recursive method which utilizes the Grasmannian framework, is robust to outliers and incomplete information but does not handle dynamic backgrounds \cite{seidel2014prost}. Since the proposed algorithm works with dynamic backgrounds, comparisons are made with  ReProCS and OPRMF only. As these algorithms are incremental in nature, they are well-equipped to handle sequential data as well as dynamic backgrounds. It will however be shown that the proposed algorithm outperforms the state-of-the-art algorithms in dynamic settings. 

Summarizing, the key contributions of the present work include (a) development of dynamic regret bounds for the inexact proximal OGD algorithm for problems with non-differentiable regularizers; and (b) development of accelerated proximal OGD algorithms in dynamic settings. The performance of all algorithms is characterized via dynamic regret bounds and via numerical tests.

		\begin{table*}
			\centering
			\captionof{table}{Dynamic Regret rates for non-smooth functions (cf. Sec.\ref{seciii})}	\label{table} 
	 \begin{adjustbox}{width=\textwidth}		\begin{tabular}{cccccc}
				\hline
				References &  Loss function & Inexact & Function class & Regret rate & Dynamic Comparator\\
				\hline 
				\cite{hall2015online} 		& Convex				& No & Non-differentiable &	$\mathcal{O}{(\sqrt{T}(1+V_{\Phi}(\boldsymbol{\theta}_T)))}$ & $V_{\Phi}(\boldsymbol{\theta}_T)):=\sum\limits_{t=1}^{T-1}\norm{\boldsymbol{\theta}_{t+1}-\Phi_t(\boldsymbol{\theta}_t)}$ \\ 
				\cite{zhang2016improved} 	& Self-concordant 	& No & \begin{tabular}{@{}c@{}}Non-smooth, \\ thrice-differentiable\end{tabular}  & 	${\mathcal{O}{(1+\min(\frac{1}{3}\mathcal{P}_T^*,4\mathcal{S}_T^*))}}$ & $\mathcal{P}_T^*=\sum\limits_{t=2}^{T}\norm{\x_t^{\star}-\x_{t-1}^\star}_t , \mathcal{S}_T^*=\sum\limits_{t=2}^{T}\norm{\x_t^{\star}-\x_{t-1}^\star}^2_t$\\
				\textbf{This work} 					& \textbf{Strongly convex}  			 	& \textbf{Yes} & \textbf{Non-differentiable}   & ${\boldsymbol{\mathcal{O}}{\boldsymbol{(1+W_K+E_K)}}}$ & $\boldsymbol{W_K:=\sum\limits_{k=2}^{K}\norm{\xks-\x_{k-1}^\star}, E_K:=\sum\limits_{k=1}^K \EE\norm{\e_k}}$\\
				\hline
			\end{tabular}	 \end{adjustbox}
			\vspace{0mm}
		\end{table*}

	\textbf{Notations:} All the scalars are denoted by small letters in regular font, vectors by  boldface small letters, and matrices by boldface capital letters.  The notation $\norm{\cdot}$ denotes standard Euclidean norm,  $\norm{\cdot}_\star$ represents a nuclear norm, and $\norm{\cdot}_F$ means the Frobenius norm. The vectorization operation on a matrix $\mathbf{A}$ is denoted by $\vect{\A}$. The Hadamard product between two matrices is denoted by the $\odot$. The all one and all zero matrices of size $m \times n$ are denoted by $\mathbf{1}_{m \times n}$ and $\mathbf{0}_{m \times n}$, respectively.

\section{Problem Formulation}\label{Prob_for}
The online convex optimization paradigm considers the online learning problem as a sequential game between a learner and an adversary \cite{zinkevich2003online}. At a given discrete time instant $k\in\N$, a learner plays an action $\x_k\in\Rn^n$. In response, an adversary selects a convex loss function $h_k:\Rn^n\rightarrow \Rn$ and the learner incurs the cost $h_k(\x_k)$. Subsequently, the adversary also reveals some information about the form of $h_k$ and the learner uses this information to determine the next action $\x_{k+1}$. Of particular interest is the scenario when $h_k$ takes the following form:
\begin{align}\label{main}
h_k = f_k + g_k,
\end{align}
where $f_k :\cX\rightarrow\Rn$ is strongly convex and differentiable with Lipschitz continuous gradient while $g_k :\cX\rightarrow\Rn $ is non-differentiable where $\cX \subset \Rn^n$ is a convex set. Such loss functions find applications in signal processing \cite{duchi2017stochastic,simonetto2015non}, robotics \cite{derenick2007convex}, and robust PCA   \cite{he2012incremental,qiu2011reprocs,kasai2016network} etc., where the smooth component corresponds to the data fitting function and the non-smooth component is the regularizer. In the context of convex optimization algorithms, such loss functions are referred to as convex composite functions and often minimized via iterative proximal point algorithms \cite{nedic2017stochastic}. In the present work however, we consider the more challenging case of time-varying loss functions and develop online learning algorithms for the same. 

When the loss function is time-varying, the learner seeks to track the minimum of $h_k$, given by $\xks = \arg\min_{\x} h_k(\x)$. Indeed, $\xks$ represents the optimal action taken by a clairvoyant that knows $h_k$ in advance. The learner however only reveals the information about  $h_k$ after taking the action $\x_k$. Specifically, the learner is revealed the full functional form of $g_k$ but only an inexact gradient $\nt f_k(\x_k):=\nabla f_k(\x_k) + \e_k$ for some $\e_k\in \Rn^n$, and uses these to find the next action $\x_{k+1}$. The performance of the learner is therefore measured by comparing the total cost incurred by the learner against that incurred by the clairvoyant, and is referred to as the dynamic regret \cite{mokhtari2016online,besbes2015non,hall2015online}
\begin{align}\label{regret}
\mathbf{Reg}_K:=\sum\limits_{k=1}^{K}[h_k(\x_k)-h_k(\xks)].
\end{align}
The notion of dynamic regret has been widely used for online learning in time-varying settings where the environment may change unpredictably over time. It may be insightful to compare $\textbf{Reg}_T$ with the related notion of static regret $\mathbf{Reg}_K^S:=\sum_{k=1}^{K}h_k(\x_k)-\min_\x\sum_k h_k(\x)$ also popular in online learning. Here, the learner's performance is measured against a static clairvoyant that knows the full sequence $\{h_k\}$ in advance. Observe that the static regret is not very useful in dynamic settings where the learner is required to be adaptive. Indeed, we always have that  $\mathbf{Reg}_K^S\leq \mathbf{Reg}_K$ and a sublinear static regret does not necessarily translate to a sublinear dynamic regret. 

In general, $\mathbf{Reg}_K$ is not sublinear in the worst case, unless restrictions are imposed on the fluctuations in $h_k$ \cite{besbes2015non} and on the gradient errors. To this end, the dynamic regret is often upper bounded in terms of regularity measures capturing cumulative variations in the optimal actions, function values, or gradient values \cite{mokhtari2016online,shahrampour2016distributed,shahrampour2018distributed}. This paper considers two such measures, namely the path length and the gradient variation, defined respectively as
\begin{align}\label{eq:var_optimal_points}
W_K&:=\sum\limits_{k=2}^{K}\norm{\xks-\x_{k-1}^\star}, \\
V_K&:=\sum\limits_{k=2}^{K}\max_{\x \in \cX}\norm{\nabla f_k(\x)-\nabla f_{k-1}(\x)}^2. \label{vkdef}
\end{align}
Intuitively, $W_K$ measures the cumulative variation of the optimal action while $V_K$ measures the cumulative gradient variation for the entire function. The two measures are related but not exactly the same, capturing complementary aspects of the problem dynamics \cite{besbes2015non,mokhtari2016online}. 

Additionally, the dynamic regret bounds also depend on the possibly stochastic error sequence $\e_k$ \cite{bedi2018tracking}. Here we develop bounds on the expected regret in terms of the cumulative mean error $E_K:=\sum_{k=1}^K \EE\norm{\e_k}$. In practical settings, the error process may arise due to various artifacts such as communication error or noise. Sec. \ref{vrsec} considers a large-scale variant of the problem, where $f_k:=\sum_{i=1}^N f_k^i$, error arises from sampling one or few of the components functions $f_k^i$ out of the $N$ components, and can be bounded in terms of $V_K$. 

The dynamic regret bounds developed here are sublinear if $W_K$ and $E_K$ are sublinear. It is therefore required that both, path variations as well as errors diminish with $k$. The path length may diminish if the target being tracked slows down over time or eventually stops. Likewise $\e_k$ may diminish if the noisy gradients available from the adversary can be corrected or improved with time. For instance, the variance reduction techniques discussed in Sec. \ref{vrsec} make the sampling error diminish with epochs. The comparison of regret result with respect to literature is summarized in Table \ref{table} where the dynamic comparators from \cite{zhang2016improved} and \cite{hall2015online} are defined respectively. The function $\Phi_t$ in \cite{hall2015online} denotes a predictor for $\boldsymbol{\theta}_t$.

 

\textbf{Robust subspace tracking:} It refers to the problem of completing a time-varying matrix that is corrupted with outliers. The static version of the problem, namely robust principal component analysis (PCA) entails decomposing a given matrix $\M$ into a sum of low-rank matrix $\L$ and a sparse matrix $\S$. The robust PCA problem has been widely studied in the context of signal processing and dimensionality reduction, and notable examples include {foreground-background separation \cite{he2012incremental}, traffic prediction \cite{qiu2011reprocs}, network anomalography \cite{kasai2016network} etc.} Interestingly, many of these tasks are inherently time-varying and the static version of the problems have always been formulated as approximations to the dynamic variants that arise in reality. For instance, the measurement matrix in the foreground-background separation problem is constructed from a window of video frames, and sliding the window over time yields the corresponding dynamic version of the problem. 

Given the measurement matrix $\M_k$ at time $k$, the cost functions for the robust subspace tracking problem are given by 
\begin{align}
f_k(\L,\S) &= \norm{\M_k-\L-\S}_F^2, \\
g_k(\L,\S) &= \lambda\norm{\L}_{\star} + \mu\norm{\vect{\S}}_1,	
\end{align}
for all $k \geq 1$, where the optimization variables $\L$ and $\S$ correspond to underlying subspace and sparse outliers, respectively. The nuclear norm regularization promotes $\L$ to have low rank structure while the $\ell_1$ norm penalty encourages $\S$ to be sparse. Additional 'box' constraints of the form $\L\in\cL$ and/or $\S \in \cS$ may be imposed on the entries of $\L$ and/or $\S$ to ensure that the values correspond to physical quantities, e.g., pixel color in videos. The goal is to track the underlying subspace and the outlier matrices,
\begin{align}
(\L_k^\star,\S_k^\star) = \arg\min_{\L,\S} f_k(\L,\S) + g_k(\L,\S), \label{rst}
\end{align}
assuming that they change slowly over time. A related problem of dynamic low-rank matrix completion was first studied from the perspective of adaptive algorithms in \cite{tripathi2017adaptive}. The analysis and algorithms presented here are however significantly more general. 

It is remarked that the classical static regret minimization framework is ill-equipped to handle such variations in the underlying subspace. Instead, static regret is more suited to the problem of online robust subspace learning, where the goal is to learn the underlying $\L$ and $\S$ matrices from sequentially available measurements. 

\section{Proposed Online Algorithm: Assumptions and Performance}\label{seciii}
For differentiable dynamic problems, online gradient descent and its variants have traditionally been the methods of choice, achieving optimal dynamic regret for many cases \cite{yang2016tracking} and incuring low complexity. While non-differentiable functions can always be handled via subgradient descent variants  \cite{ISG_1,ISG_2,ISG_3}, it is known from the convex optimization literature, that convex composite functions (differentiable+non-differentiable) are better handled via proximal point methods  \cite{bertsekas2011incremental,nedic2017stochastic}. We build upon this intuition and put forth an inexact proximal OGD (IP-OGD) algorithm that takes the form:
\begin{align}\label{prox_update1}
\x_{k+1} = \pk{\x_{k} - \alpha \nt f_{k}(\x_{k})}, 
\end{align}
for all $k\geq 1$. The proximal operator is defined with respect to non-differentiable component $g_k(\cdot)$ as
\begin{align}\label{proximal}
\pk{\x} = \arg\min_{\u\in\cX}  g_{k}(\u) + \frac{1}{2\alpha} \norm{\u-\x}^2_2 ,
\end{align}
where $\alpha>0$ is the step-size parameter. The full algorithm is summarized in Algorithm \ref{algo_1}.  
	  \begin{algorithm}
	  	\caption{IP-OGD Algorithm}\label{algo_1}
	  	\begin{algorithmic}[1]
	  		\STATE {\textbf{Initialize}} $\x_{1}$	  		
	  		\STATE {\textbf{for} k = 1, 2,$\cdots$, $K$} \textbf{do} 	  		
	  		\STATE \ \ \ \ \textbf{Action} $\x_k$  	  		
	  		\STATE \ \ \ \ \textbf{Observe} inexact gradient $\nt f_k(\x_k)$ and function $g_k$  		
	  		\STATE \ \ \ \ \textbf{Update} $\x_{k+1} = \pk{\x_{k} - \alpha \nt f_k(\x_k)}$   		
	  		\STATE \textbf{end for }
	  	\end{algorithmic}
	  		  	\label{algo:IP-OGD}
	  \end{algorithm}
Having stated the algorithm, we explicate the key assumptions required for establishing the regret bounds. 

\subsection{Assumptions}
The development of the regret bounds relies on the following key assumptions: 
	  \begin{itemize}
	  	\item[\textbf{A1.}] \textbf{Lipschitz continuity:} the function $f_k$ is $L$-smooth on $\cX$, i.e., 
	  	\begin{align}\label{Lip_f1}
	  	\norm{\nabla f_k(\x) - \nabla f_k(\y)}\leq L\norm{\x-\y} , 
	  	\end{align}
	  	for all $k\in\N$ and $(\x,\y)\in\cX$. The function $g_k$ is Lipschitz continuous on $\cX$  with parameter $L_g>0$, i.e., 
        	\begin{align}\label{Lip_g1}
        \norm{ g_k(\x) -  g_k(\y)}\leq L_g\norm{\x-\y}.  
        \end{align}
        for all $k\in\N$ and $(\x,\y)\in\cX$. 
	  	\item[\textbf{A2.}] \textbf{Strong convexity:}  The function $f_k$ is $\mu$-convex, i.e.,
	  		  	\begin{align}\label{eq:SC}
	  	    \ip{\nabla f_k(\x) - \nabla f_k(\y),\x-\y} \geq \mu \norm{\x-\y}^2_2  , 
	  	\end{align}
for all $k\in\N$ and $(\x,\y)\in\cX$.			
	 \item[\textbf{A3.}] \textbf{Bounded variations:}  the distance between any two consecutive optimal points is upper bounded by the following condition 
	 \begin{align}\label{eq:BV}
	 \norm{\xks- \x_{k-1}^\star} \leq \sigma ,
	 \end{align}
	for all $k$. Also, the error process has bound mean squared error, i.e., 
	 \begin{align}\label{eq:BV2}
	 \EE\norm{\e_k} \leq \gamma. 
	 \end{align}
\end{itemize}  

It can be seen that the assumptions are quite standard and are satisfied for most problems of interest. Specifically, the smoothness and strong convexity properties of $\nabla f_k$ are necessary for developing a contraction relation between $\norm{\x_{k+1}-\x_k^*}$ and $\norm{\x_k-\x_k^*}$. The Lipschitz continuity of $g_k$ gives a bound on its sub-gradient which will be useful in bounding the sub-gradient of the combined convex function $h_k$. The strong convexity of $f_k$ also makes the combined function $h_k$ strongly convex, ensuring that the optimum $\xks$ is unique. Strongly convex objectives arise in a number of machine learning applications, such as robust subspace tracking, Lasso, support vector machines, etc. Even for other applications where the objective function is not strongly convex, a regularization of the from $\frac{\mu}{2}\norm{\x}^2_2$ may at times be included without sacrificing the performance significantly. The assumption on bounded variations of the optimum $\xks$ are natural for most tracking applications \cite{derenick2009convex,derenick2009optimal}. For instance, in the context of robust subspace tracking, \textbf{(A3)} ensures that the underlying target subspace does not change abruptly. 

						
						
\subsection{Regret Bounds}\label{vr} 
We now present the regret bounds for the IP-OGD algorithm. The development of the regret bounds entails three key technical results, namely an upper bound on the post-update optimality gap $\norm{\x_{k+1}-\xks}$ in terms of the pre-update gap ${\norm{\x_{k}-\xks}}$, an upper bound on the cumulative optimality gap $\sum_k \norm{\x_k-\xks}$, and an upper bound on the subgradient norm $\partial h_k(\x_k)$. The first two bounds are presented directly in terms of $\norm{\e_k}$ while the bound on the subgradient norm and the final regret bound are provided in expectation. 

The first lemma is crucial and quantifies the impact of a single proximal update. At iteration $k$, the subsequent lemma bounds $\norm{\x_{k+1}-\xks}$ in terms of $\norm{\x_k-\xks}$ and the gradient error $\e_k$.   
	  		\begin{lem}\label{lem1}
	  			Under \textbf{(A1)-(A2)}, the sequence of  $\{\x_k\}_{k\in\N}$ generated by the IP-OGD algorithm satisfies 
	  			\begin{align}\label{Lemma_1}
	  			&\norm{\x_{k+1} - \xks}\leq \rho \norm{\x_{k} - \xks}+\alpha\norm{\e_k},
	  			\end{align}
	  				where, $\rho^2:= \left(1  - \alpha(2\mu - \alpha L^{2})\right)$.
	  		\end{lem}
The proof of Lemma \ref{lem1} is provided in Appendix \ref{proof:lemma1}. The bound is useful when $\rho < 1$ or equivalently when $\alpha \in (0, 2\mu/L^2)$. For such a choice of $\alpha$, the subsequent lemma develops a bound on the cumulative optimality gap. 
\begin{lem}\label{lem2}
	Under \textbf{(A1)-(A2)} and for $0<\alpha < 2\mu/L^2$, the sequence of $\{{\x}_k\}_{k\in\N}$ generated by the IP-OGD algorithm satisfies
	  			\begin{align}
	  			\!\!\!\!\!\!\sum\limits_{k=1}^{K}\norm{{\x}_{k}- \xks }  \leq & \frac{1}{1-\rho}\left[\!\norm{{\x}_{1} \!- \x_{1}^\star} \!+\! W_K  \!+\!  \alpha\sum\limits_{k=1}^{K} \norm{\e_k}\right]. 
	  			\end{align}
\end{lem}	
The proof of Lemma \ref{lem2} is provided in Appendix \ref{proof:lemma1}, and uses the  inequality $\norm{\x_k-\xks} \leq \norm{\x_k-\x_{k-1}^\star} + \norm{\xks-\x_{k-1}^\star}$ and the result of Lemma \ref{lem1} as the key ingredients. Finally the subgradient norm can be bounded through the use of \textbf{(A3)} in the following lemma whose proof is provided in Appendix \ref{proof:lemma1}. 
\begin{lem}\label{lemsg} 
Under \textbf{(A1)-(A3)} and for $0<\alpha < 2\mu/L^2$, the subgradient norm of the objective function is bounded as
\begin{align}
\!\!\!\!\EE\norm{\partial h_k(\!\x_k\!)}\!\leq\! D\!:=\!\frac{1\!+\!\rho}{\alpha} \bigg(\!\!\!\norm{{\x}_{1} \!\!- \!\x_{1}^\star} \!+\!  \frac{\alpha \gamma \!+\! \sigma}{1\!-\!\rho} \bigg) \!\! +\!\! 2\gamma \!\!+\!\! 2L_g,
\end{align}
for all $k\in\N$. 
\end{lem} 

Finally, the following theorem provides the bound on the expected regret. 
\begin{thm}\label{thm1}
Under \textbf{(A1-A3)} and for $0<\alpha < 2\mu/L^2$, the dynamic regret of the IP-OGD algorithm is bounded as
\begin{align} \label{regret1}
	\Ex{\mathbf{Reg}_K}\leq \mathcal{O}(1+W_K+E_K).
\end{align}
\end{thm}

 \begin{IEEEproof}
From first order convexity condition, we have 
	\begin{align}
 	\sum\limits_{k=1}^{K}[&h_k(\x_k)-h_k(\xks)]\leq  \sum\limits_{k=1}^{K}\ip{\partial h_k(\x_k), \x_k-\xks}\label{proof:lem2:first}
	  		 	\\
	  		 	&\leq \sum\limits_{k=1}^{K}\norm{\partial h_k(\x_{k})}  \norm{\x_k- \xks}\label{proof:lem2:second}
	  		 	\\
	  			&\leq \sum\limits_{k=1}^{K}D \norm{\x_k- \xks}\label{proof:lem2:third}
	  		 		\\
	  		 	&\leq  \frac{D}{1\!-\!\rho}\!\norm{\x_{1} \!-\! \x_{1}^\star} \!+\! \frac{D}{1-\rho}W_{K} + \frac{\alpha D}{1-\rho}\sum_{k=1}^{K} \norm{\e_k}.\label{last} 		 	
	  		 	\end{align}
The inequality in \eqref{proof:lem2:second} follows form the Cauchy- Schwartz inequality. Next inequality in \eqref{proof:lem2:second} follows from the gradient boundedness as proved in Proposition \ref{lemsg}. The last inequality in \eqref{last} follows from the result in Lemma \ref{lem2}. Taking expectations on both sides of \eqref{last}, we obtain the required result. \end{IEEEproof}

The regret obtained in Theorem \ref{thm1} is sublinear as long as both $W_K$ and $E_K$ are sublinear. Further, if the gradient error is zero, we obtain the dynamic regret of $\mathcal{O}(1+W_K)$ which is of the same rate as that obtained in \cite{mokhtari2016online} for strongly convex differentiable functions. In case the path length is not sublinear, the corresponding tracking performance can be readily obtained {as provided in Lemma \ref{tracking_performance}}. It is remarked that unlike in  \cite{shahrampour2018distributed} and \cite{pmlr-v38-jadbabaie15}, the step size parameter can be selected within a specified range and need not depend on the  $W_K$. 

\begin{lem}\label{tracking_performance}
	Under \textbf{(A1)-(A3)}, the sequence of iterates $\{\x_k\}$ generated by IP-OGD algorithms satisfies 
	\begin{align}
	\limsup_{k\rightarrow\infty} 	\EE\norm{\x_k - \xks} \leq   \frac{\alpha \gamma + \sigma}{1-\rho}. \label{xkbndf1}
	\end{align}
\end{lem} 
This result follows from \eqref{xkbndf0} after taking limit $k\rightarrow\infty$, This result states that the proposed algorithms tracks the adversarial target with some error. 


\section{IP-OGD for Large Scale Learning}\label{vrsec} 
The inexact algorithm proposed in Sec. \ref{seciii} opens up new avenues for development of variance reduced gradient descent methods in dynamic settings. These algorithms are useful in large-scale problems when the loss function $f_k$ is expressible as a sum of several component functions, 
\begin{align}\label{large_Scale}
 f_k(\x):=\frac{1}{N}\sum_{i=1}^{N}f_k^i(\x).
\end{align} 
For such loss functions, it may not be possible to calculate the full gradient $\nabla f_k$ at every iteration. For instance, in the robust subspace tracking problem \eqref{rst}, each component function may correspond to a frame and processing a large number of frames at every iteration may not be viable. Variance reduction techniques rely on calculating the inexact gradient $\nt f_k$ that uses only one component function $f_k^i$ at every iteration. A correction factor is added to the inexact gradient that utilizes an old value of the full gradient, making the approximation better with iterations. Variance reduction has been widely applied to proximal gradient methods for solving static convex optimization problems \cite{xiao2014proximal,meng2017asynchronous}. 

However, achieving variance reduction in dynamic settings is challenging and has never been attempted before. To illustrate the challenge, consider a setting where $\nabla f_k^1$ is utilized at iteration $k > 1$. Due to the time-varying nature of the function, no information about the other components $\{f_k^i\}_{i=2}^N$ can be used at the current or subsequent iterations. In other words, sampling the component functions leads to loss of information. We alleviate the difficulty by exploiting the fact that gradient variations diminish over iterations. Consequently, progressive variance reduction is achieved from using gradients of stale component functions. 

We put forth two large-scale learning algorithms inspired from SVRG  \cite{xiao2014proximal} and ISSG \cite{friedlander2012hybrid}. The algorithms are motivated as special cases of the IP-OGD algorithm with the gradient error arising from the sampling and correction process specific to each algorithm. Bounds on the cumulative error $E_K$ are developed for both cases, and the overall regret bounds follow simply from plugging in the corresponding bounds for $E_K$ into \eqref{regret1}. 

\subsection{Assumptions}
We begin with discussing the assumptions required for the proposed large-scale learning algorithms. Different from the IP-OGD setting detailed earlier, the component gradient functions $\nabla f_k^i$ are no longer inexact. However the gradient used at every iteration is still an approximate version of the full gradient with error depending the algorithms used. Given the finite sum structure of the loss function, bounds on $E_K$ are developed in terms of the path length $W_K$ and the cumulative gradient variation of each component $f_k^i$, defined as
\begin{align}\label{vki}
V_K^i&:=\sum\limits_{k=2}^{K}\max_{\x \in \cX}\norm{\nabla f_k^i(\x)-\nabla f_{k-1}^i(\x)}^2.
\end{align}
It is remarked that the standard gradient variation $V_K$ defined in \eqref{vkdef} for the function $f_k$ does not suffice since it does not capture variations in individual $f_k^i$s. Indeed, $V_K$ may be zero even when $\{V_K^i\}$ are non-zero. Further, the assumptions in \textbf{A1}-\textbf{A3} are modified as follows: 
\begin{itemize}
	\item[\textbf{A1'.}] \textbf{Lipschitz continuity.} The functions $\{f_k^i\}$ are all $L$-smooth on $\cX$, i.e., $\norm{\nabla f_k^i(\x)-\nabla f_k^i(\y)} \leq L\norm{\x-\y}$ for all $\x,\y \in \cX$.  
	\item[\textbf{A2'.}] \textbf{Strong convexity.} The functions $\{f_k^i\}$ are all  $\mu$-convex on $\cX$, i.e., $\ip{\nabla f_k^i(\x)-\nabla f_k^i(\y),\x-\y} \geq \mu\norm{\x-\y}^2$ for all $\x,\y \in \cX$.
	\item[\textbf{A3'.}] \textbf{Bounded gradients.} The functions $\{f_k^i\}$ have bounded gradients, i.e., 
	\begin{align}\label{eq:BVg}
	\max_{i,k}\norm{\nabla f_{k}^{i}(\x_k)} \leq M. 
	\end{align}
\end{itemize}
Assumptions \textbf{(A1')-(A2')} are similar to \textbf{(A1')-(A2')}, but applied to each component separately. Assumption \textbf{(A3')} is however different, though widely used in convex optimization. It is remarked that \textbf{(A3')} holds if $\cX$ is a compact. To see this, observe that since $f_k^i$ is convex, it is also locally Lipschitz, and consequently Lipschitz continuous on the compact set $\cX$, implying that $\norm{\nabla f_k^i}$ is bounded. Conversely, if $\cX$ is not bounded and there exist some $\x, \y \in \cX$ such that $\norm{\x-\y} \rightarrow \infty$, the strong convexity assumption \textbf{(A2')} on $f_k^i$ would preclude the possibility of the gradient being bounded.  


We begin with discussing the online SVRG algorithm, that relies on old gradients to construct an approximation for the full gradient at every iteration, while using only a single component $\nabla f_k^i$. As in the classical variance reduced algorithms, the proposed online algorithms also require additional memory. The two algorithms differ in the manner in which stored gradients are updated and consequently offer complementary advantages in terms of memory usage, computational complexity, and dynamic regret bounds.

\subsection{Online proximal-SVRG algorithm}\label{secsvrg}
The OP-SVRG is based on the proximal SVRG algorithm from \cite{defazio2014saga} and entails storing $\m$, the mean value of obsolete gradients. In order to ensure that the gradients are not too old, $\m$ is recalculated every few iterations. The algorithm starts by evaluating $\bt^i = f_1^i(\x_1)$ for all $1\leq i \leq N$ and some $\x_1 \in \cX$. Further, $\m:=\frac{1}{N}\sum_j \bt^j$ is calculated and stored in the memory, along with the current iterate $\tilde{\x}=\x_1$. Subsequently at time $k$, both
$\nabla f_k^i(\x_k)$ and $\nabla f_k^i(\tilde{\x})$ are evaluated and the IP-OGD updates are carried out with the inexact gradient:
\begin{align}
\nt f_k(\x_k) = \nabla f_k^i(\x_k) - \nabla f_k^i(\tilde{\x}) + \m,
\end{align}
where the component $i$ is selected either randomly or in a cyclic manner. After every $K_0$ iterations, $\tilde{\x}$ is refreshed, all gradients recalculated, and $\m$ is updated. The full algorithm implementation is summarized in Algorithm \ref{algo:SVRG}. 

   \begin{algorithm}
	  	\caption{Online Proximal-SVRG}\label{algo:SVRG}
	  	\begin{algorithmic}[1]
	  		\STATE {\textbf{Initialize}} $\x_1$ 
	  		\FORALL{$k = 1,\ldots, K$}  		
				\IF{$k (\text{mod } K_0) = 1$}
					\STATE \textbf{Update} $\tilde{\x} \leftarrow \x_k$
					\STATE \textbf{Set} $\m = \frac{1}{N}\sum_{j=1}^N\nabla f_k^j(\tilde{\x})$ 
				\ENDIF 
	  		\STATE \textbf{Select} $i \in \{1, \ldots, N\}$ \label{isvrg}  	  		
	  		\STATE \textbf{Observe} exact gradients $\nabla f_k^i(\x_k)$ and $\nabla f_k^i(\tilde{\x})$ 
	  	   \STATE \textbf{Update} $\z_{k} = \x_{k} - \alpha \big(\nabla f_k^{i}(\x_k)-\nabla f_k^i(\tilde{\x}) + \m \big) $   		
	  		\STATE \textbf{Update} $\x_{k+1} = \pk{\z_k}$
				\ENDFOR
				\end{algorithmic}
	\end{algorithm}	  
	
The gradient error for OP-SVRG algorithm is given by 
\begin{align}
\!\!\!\e_k \!= \!\nabla f_k^i(\x_k)\!-\!\nabla f_k^i(\x_{\tau_k}) \!\!+\!\! \frac{1}{N}\!\!\sum\limits_{j=1}^{N}(f_k^j(\x_{\tau_k}) \!-\! \nabla f^j_k(\x_k)), \label{error_svrg}
\end{align}
where $\tau_k = \lfloor \frac{k}{K_0}\rfloor K_0 + 1$ is the value of $k$ corresponding to the last memory update. Again, the index $i$ may be selected in either random or cyclic fashion. Since the OP-SVRG is again motivated as an inexact algorithm, its analysis culminates in the following bound on error. 

\begin{lem}\label{lemsvrg}
For $\mu/L > 0.89$, there exists $\alpha < 2\mu/L^2$ such that the OP-SVRG gradient error satisfies
\begin{align}
E_K = \mathcal{O}\left(1+W_K+\sum\limits_{j=1}^{N}\sqrt{KV_K^j}\right). \label{eksvrg}
\end{align}
\end{lem}
The expected regret bound for OP-SVRG is obtained from plugging the bound in \eqref{eksvrg} into the expected regret bound for the inexact algorithm in \eqref{regret1} and yields $\Ex{\mathbf{Reg}_K}\leq \mathcal{O}(1+W_K+\sum\limits_{j=1}^{N}\sqrt{KV_K^j})$. It is remarked that the result in \eqref{eksvrg} requires the problem to be well-conditioned with $\mu > 0.89L$.

The online proximal SVRG offers significant computational advantages over IP-OGD. To see this, let us assume that the calculation of $\nabla f_k^i$ and application of the proximal operator to an $n$-dimensional vector incur $\mathcal{O}(n)$ complexity. Then the complexity incurred by Algorithm \ref{algo:SVRG} is $\mathcal{O}(Kn(2+\frac{N}{K_0}))$, in contrast to the $\mathcal{O}(KNn)$ complexity of IP-OGD. Clearly, the cost savings from SVRG are significant only when $K_0 \gg 1$. As compared to IP-OGD, Algorithm \ref{algo:SVRG} also requires additional memory for storing the vector $\m$. 

\subsection{Online Proximal gradient method with increasing sample size}
The dynamic regret rates obtained in \ref{secsvrg} can be improved at the cost of additional computational complexity by making use of an incremental update approach with increasing sample size. In contrast to the classical incremental algorithms that entail carrying out the updates using only a single component $\nabla f_k^i$, the idea is to use increasing number of component functions at every iteration. Specifically, the inexact gradient used at iteration $k$ takes the form:
\begin{align}
\nt f_k(\x_k) = \frac{1}{\abs{\I_k}}\sum_{j \in \I_k} f_k^j(\x_k),
\end{align}
where $\I_k \subset \{1, \ldots, N\}$ is selected randomly. The full algorithm is summarized in Algorithm \ref{algoiss}. 
\begin{algorithm}
	  	\caption{Online Proximal with Incremental Sampling Size}\label{algoiss}
	  	\begin{algorithmic}[1]
	  		\STATE {\textbf{Initialize}}  $\x_1$	  	
				\FORALL{$k = 1, \ldots, K$}
				\STATE \textbf{Select} $\I_k \subset \{1, \ldots, N\}$	
				\STATE \textbf{Calculate} $\nt f_k(\x_k) := \frac{1}{\abs{I_k}}\sum\limits_{i \in I_k}\nabla f_k^i(\x_k) $ 
	  	   \STATE \textbf{Update} $\z_{k} = \x_{k} - \alpha\nt f_k(\x_k)$   		
	  		\STATE \textbf{Update} $\x_{k+1} = \pk{\z_{k}}$
	  		\ENDFOR
	  	\end{algorithmic}
	  \end{algorithm}
  
The sampling process in Algorithm \ref{algoiss} ensures that the mean square of the gradient error given by 
 \begin{align}
  \e_k=\frac{1}{\abs{\I_k}}\sum\limits_{i \in \I_k}\nabla f_k^i(\x_k)- \frac{1}{N}\sum\limits_{i=1}^{N}\nabla f_k^{i}(\x_k), \label{error_incresing}
  \end{align} 
diminishes with $k$. The following proposition provides a choice of samples sizes that lead to an $\mathcal{O}(1)$ bound on $E_K$. 
	
 \begin{lem}\label{lemiss}
For the choice $\abs{\I_k} = \lceil N(1-e^{-k}) \rceil$ in Algorithm \ref{algoiss}, the bound on the cumulative error takes the form 
\begin{align}
E_K \leq \mathcal{O}(1).
\end{align} 
 \end{lem}

The proof of Lemma \ref{lemiss} is provided in Appendix \ref{proof_prop_3}, and is applicable for any choice of $\abs{\I_k}$ such that $N-\abs{\I_k}$ is summable. Plugging back the result of Lemma \ref{lemiss}, we obtain the expected regret of Algorithm \ref{algoiss} as $\mathcal{O}(1+W_K)$. Note that achieving constant $E_K$ requires the evaluation of $\mathcal{O}(\abs{\I_k})$ gradients per iteration, leading to an overall complexity of $\mathcal{O}(NKn)$.

\section{Numerical Results}
This section provides numerical performance results of the proposed algorithms. The performance of the IP-OGD algorithm is tested on a synthetic dataset corresponding to a bot formation control problem, and the results in Theorem 1 are verified. Next, the computational advantages of the OP-SVRG algorithm are demonstrated on a robust subspace tracking problem using video dataset. 

\subsection{Coordinating bot formations}
In this problem, it is required for a group of robots to track and follow a leader bot estimating the leader's trajectory from the received signal strength at each discrete time instant. The follower bots try to maintain a given canonical pose or shape at each time instant while following the leader such that the leader bot is also a part of the pose. Such problems can generally be formulated as convex optimization problems where the goal is to minimize a cost function while satisfying certain constraints \cite{derenick2007convex}. As the bots must move instantaneously in response to the movement of the leader, solving a large convex optimization problem at every time instant is not viable. The problem is particularly challenging in dynamic settings, e.g., if the pose or shape requirement is also time-varying.

\subsubsection{Problem Formulation}
We follow the formulation in \cite{derenick2007convex} where a fusion center gathers all the data from the follower bots, takes decisions, and broadcasts them to each bot. Different from \cite{derenick2007convex}, we consider an online and time-varying setting, with problem data arriving sequentially over time. Consider a formation of $m+1$ robots in a two dimensional plane and let $\p_i(k) \in \Rn^2$ be their position at time $k$. Without loss of generality, the leader is indexed as 1, and all position vectors are concatenated into the matrix $\p(k) \in \Rn^{(m+1) \times 2}$. The desired shape of the formation at time $k+1$ relative to the leader's position at $\p_1(k)$ is described by a set of linear equations of the form
\begin{align}\label{const}
\A(k) \p(k+1) = \textbf{0}.
\end{align}
where $\A(k) \in \Rn^{(2m-2) \times (2m+2)}$ and $\textbf{0}$ is the all zero matrix of appropriate size. At time $k$,  robots estimate the position of the leader $\p_1(k+1)$ and follow it while maintaining the shape described by \eqref{const}. {The further details about the constraint in \eqref{const} can be found in \cite{derenick2007convex} and discussed in  Appendix \ref{bot_form} of the supplementary material.} To this end, the leader broadcasts a signal that is subsequently used by the robots in order to estimate $\p_{1}(k+1)$. Denoting the measurement received by robot $i$ at time $k$ by $z_i(k)$, the robots attempt to solve the following constrained weighted least squares problem in order to determine $\p(k+1)$:
\begin{align}
\p(k+1):=\arg\min_{\p} &\sum\limits_{i=2}^{m+1} \frac{(z_i(k) - \v_i^T\p_1)^2}{\sigma_i^2}  + \lambda \norm{\p-\p(k)}_2^2 \nonumber\\
\text{s. t.} & \hspace{1cm}\A(k)\p = \textbf{0} \label{botform}
\end{align}
where $\v_i$ denotes the regressor corresponding to the measurement $\z_i$ and $\sigma_i^2$ denotes the corresponding noise variance estimate. A squared movement penalty term is added to minimize unnecessary movements and encourage faster response times. 

In order to apply the IP-OGD method, we observe that
\begin{align}
f_k(\p) &= \sum\limits_{i=2}^{m+1} \frac{(z_i(k) - \v_i^T\p_1)^2}{\sigma_i^2} + \lambda \norm{\p-\p(k)}^2,  \\
g_k(\p) &= \ind_{A(k)\p = 0}(\p),
\end{align}
where the indicator function is 1 when $\A(k)\p = 0$ and infinity otherwise. Consequently, the proximal OGD updates take the form:
\begin{align}
\y_k &= \p(k) - \alpha \nabla f_k(\p(k)) \\
\p(k+1) &= \prox_{g_k}^\alpha(\y_k) \\
&= (\mathbf{I}-\A(k)^{\dagger}\A(k))\y_k,
\end{align}
where $\A(k)^{\dagger}$ denotes the Moore-Penrose pseudo inverse of $\A(k)$.

\begin{figure}[h]
	\centering
	\includegraphics[scale=0.4]
	{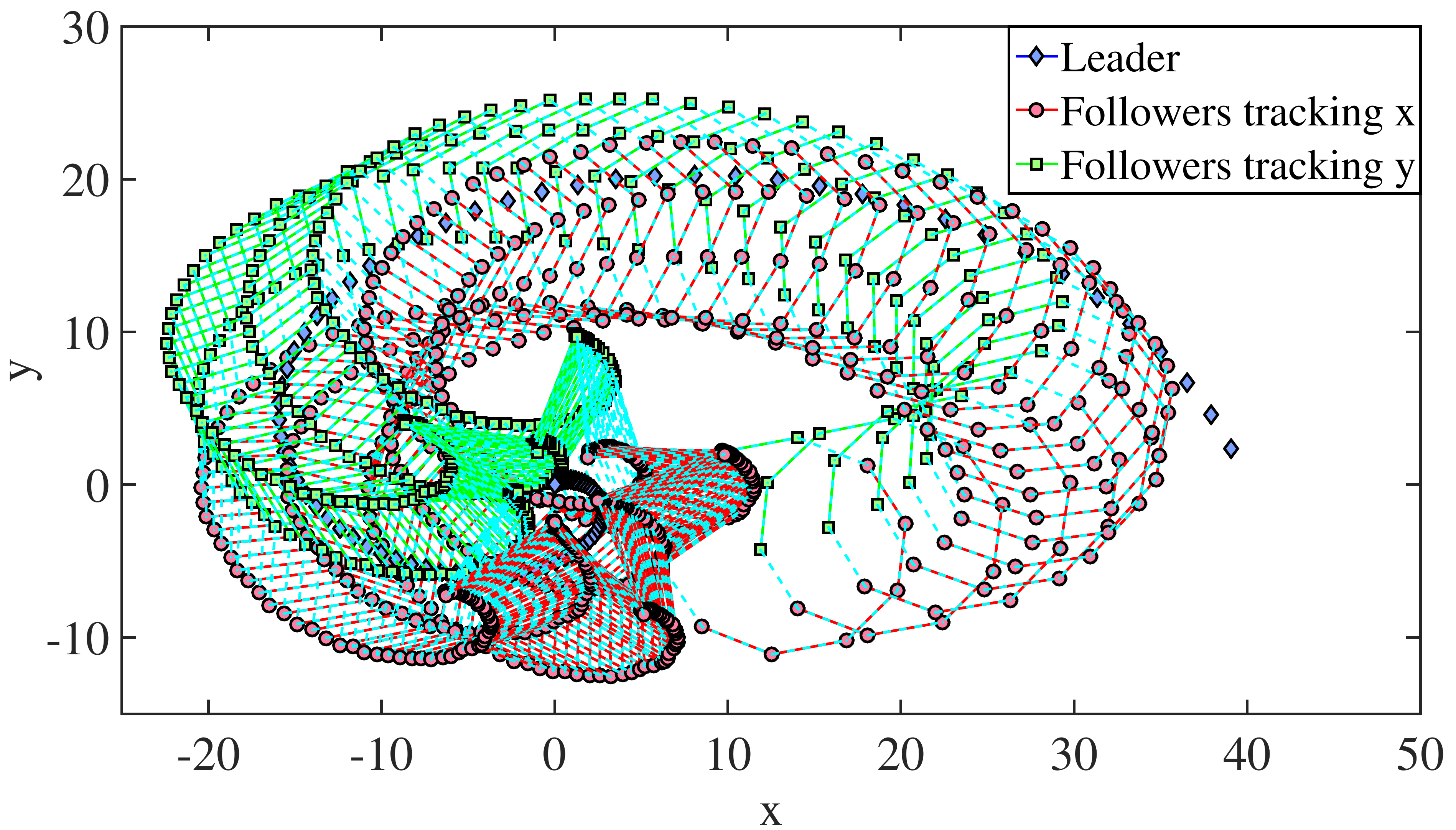}
	\caption{Example run of the IP-OGD algorithm for $m = 10$ robots. The shape of the formation morphs from a decagon to a star.}
	\label{figstar}
\end{figure}

\begin{figure}[h]
	\centering
		\includegraphics[scale=0.4]
	{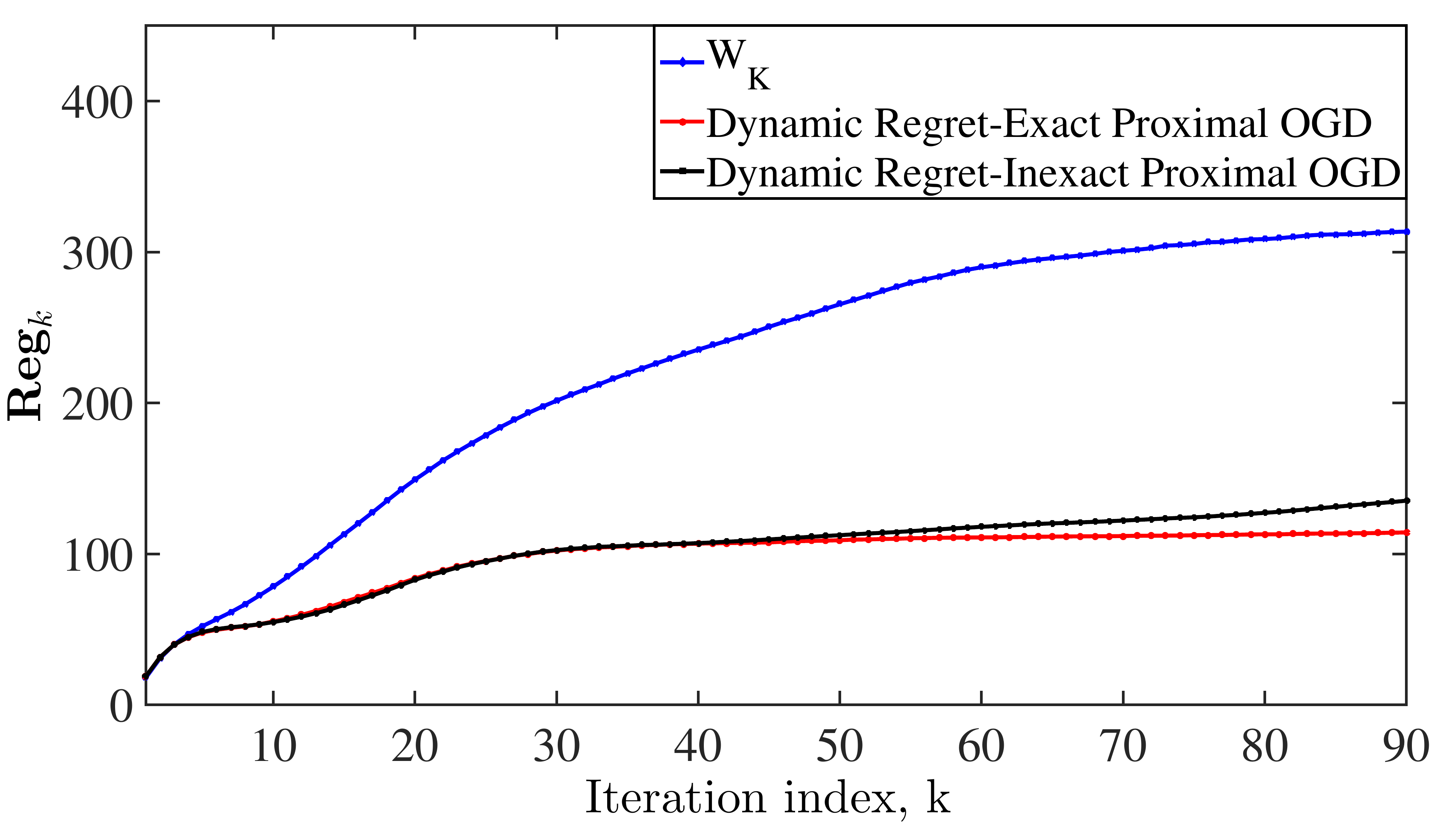}
	\caption{Dynamic regret for IP-OGD algorithm}
	\label{figreg}
\end{figure}

\subsubsection{Simulation Results} 
For the simulations, the leader bot trajectory $\p_1(k)$ is varied according to the following spiral dynamics $\p_1(k) = (1-k/K)^2[40\cos(\omega k) ~ 30\sin(\omega k)]^T$, where $\omega = 0.06$ and $K$ is the total number of times for which tracking is performed. A total of $m=10$ robots track the leader, of which 6 bots track the $x$-coordinate ($\v_i = [1~0]^T$ for $1\leq i\leq 6$) while the other 4 track the y-coordinate ($\v_i = [0~1]^T$ for $7\leq i\leq 10$). The received measurements are modeled as $z_i(k) = \v_i^T\p_1(k+1) + n_i$, where each $n_i \sim \mathcal{N}(0,\sigma_i^2)$ with $\sigma_i^2 = 0.01$. We also set $\lambda = 3$ and $\alpha = 0.2$. The shape starts as a regular decagon and gradually changes into a non-convex star shape. 


An example run of the algorithm is shown in Fig. \ref{figstar}. The dynamic regret for this case is plotted in Fig. \ref{figreg}. It can be seen that the path length as well as the dynamic regret for IP-OGD algorithm is sublinear. Finally, Gaussian noise with variance 0.3 is added to the exact gradient.


%

\subsection{Dynamic foreground background separation}
We consider the robust subspace tracking problem detailed in Sec. \ref{Prob_for} and apply it for foreground-background separation in video. At any given time a total of $L$ video frames $\{\m_k \in \Rn^r\}_{k=1}^L$ are observed and collected into a matrix $\M_k$. The goal is to split the matrix into low-rank component $\L_k$ and a sparse component $\S_k$ at each time $k$. Towards this end, we consider the following loss function and regularizer:
\begin{align}
f_k(\L,\S) &= \norm{\M_k-\L-\S}_F^2 + \mu_{\L}\norm{\L}_F^2 + \mu_{\S}\norm{\S}_F^2, \\
g_k(\L,\S) &= \lambda_{\L}\norm{\L}_{\star} + \lambda_{\S}\norm{\vect{\S}}_1 ,
\end{align}
Additional ridge regularizers are added to $f_k$ in order to improve its condition number.

We first consider the IP-OGD algorithm without gradient errors. The updates for IP-OGD take the form:
\begin{align}
\mathbf{Z}_{k+1}&= \mathbf{L}_k- \alpha_{\mathbf{L}}\bigg( 2 (\mathbf{L}_k+\mathbf{S}_k-\mathbf{M}_k) +2 \mu_{\mathbf{L}} \mathbf{L}_k\bigg)\nonumber \\
\mathbf{L}_{k+1}&=\D_{\alpha_{\mathbf{L}}\lambda_{\mathbf{L}}}(\mathbf{Z}_{k+1}), 
\end{align}
\begin{align}
\mathbf{Y}_{k+1}&= \mathbf{S}_k- \alpha_{\mathbf{S}}\bigg( 2 (\mathbf{L}_{k}+\mathbf{S}_k-\mathbf{M}_k) +2 \mu_{\mathbf{S}} \mathbf{S}_k\bigg) \nonumber \\
\mathbf{S}_{k+1}&= \mathcal{S}_{\alpha_{\mathbf{S}} \lambda_{\mathbf{S}}}(\mathbf{Y}_{k+1}) ,
\end{align}
where the singular value thresholding and shrinkage operators are defined as follows. For a matrix $\Z = \U\Sig\V^T$, we have that $\D_\lambda(\Z) = \U\D_\lambda(\Sig)\V^T$ where $\D_\lambda(\Sig)$ is a diagonal matrix whose entries are given by $[\D_\lambda(\Sig)]_{ii} = \max\{[\Sig]_{ii}-\lambda,0\}$. Likewise, the entries of $\D_\lambda(\Y)$ are given by 
\begin{align}
[\D_\lambda(\Y)]_{ij} &= \begin{cases} [\Y]_{ij}-\lambda, & [\Y]_{ij} > \lambda \\
0, & -\lambda < [\Y]_{ij} < \lambda \\
 [\Y]_{ij} + \lambda, & [\Y]_{ij} < -\lambda .
\end{cases}
\end{align}
For a matrix $\Y$, the shrinkage operator $\mathcal{S}_{\lambda}(\Y)$ is a matrix whose entries are given by
\begin{align}
[\mathcal{S}_{\lambda}(\Y)]_{ij} = \text{sign}([\Y]_{ij})(\abs{[\Y]_{ij}}-\lambda)_{+}
\end{align}

\begin{figure}
\centering
	\includegraphics[scale=0.6]{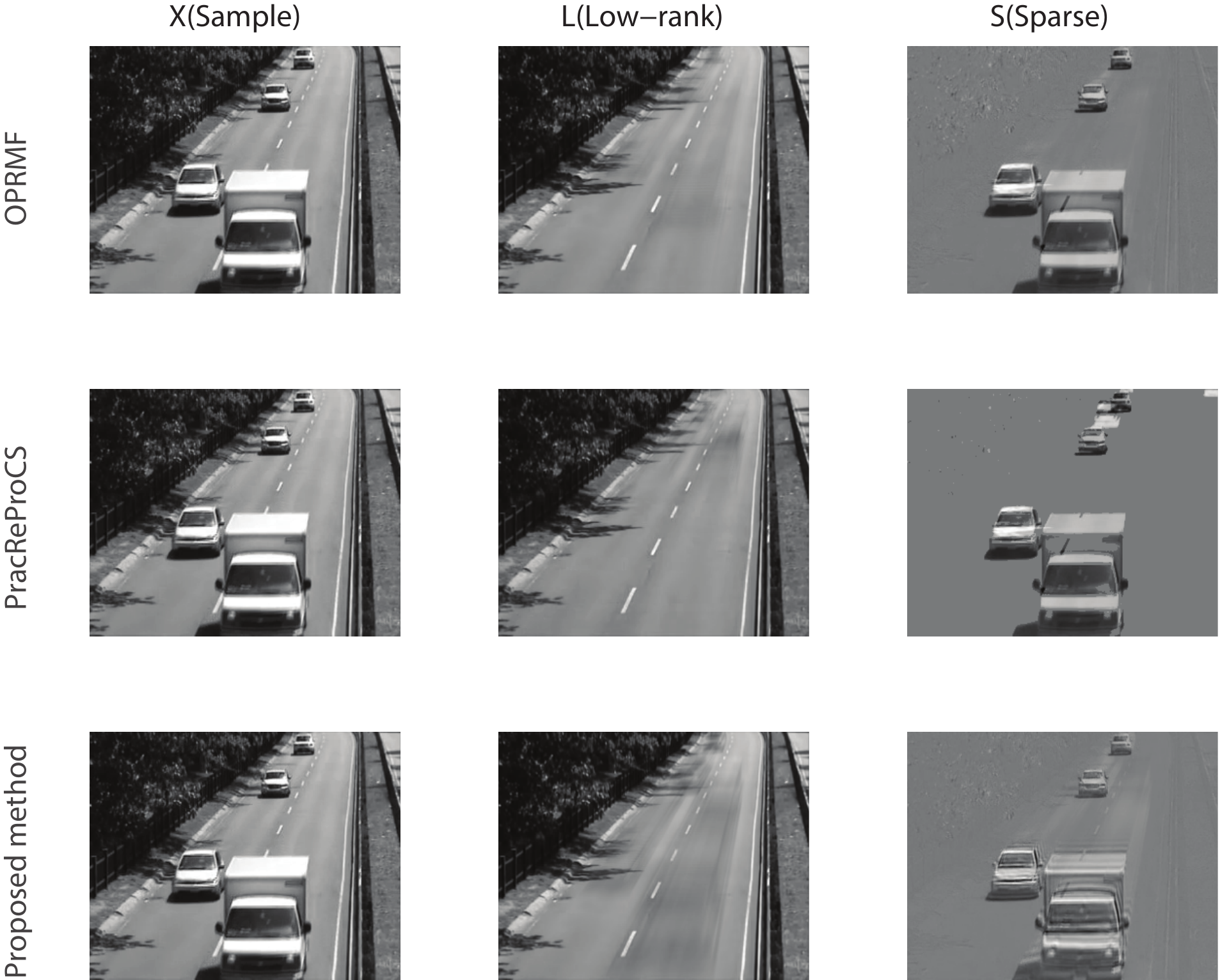}
	\caption{Performance comparisons with OPRMF \cite{wang2012probabilistic} and ReProCS \cite{reprocs} over static background.}\label{static}
\end{figure}
\begin{figure}
\centering
\hspace*{-0.5cm}	\includegraphics[scale=0.6]{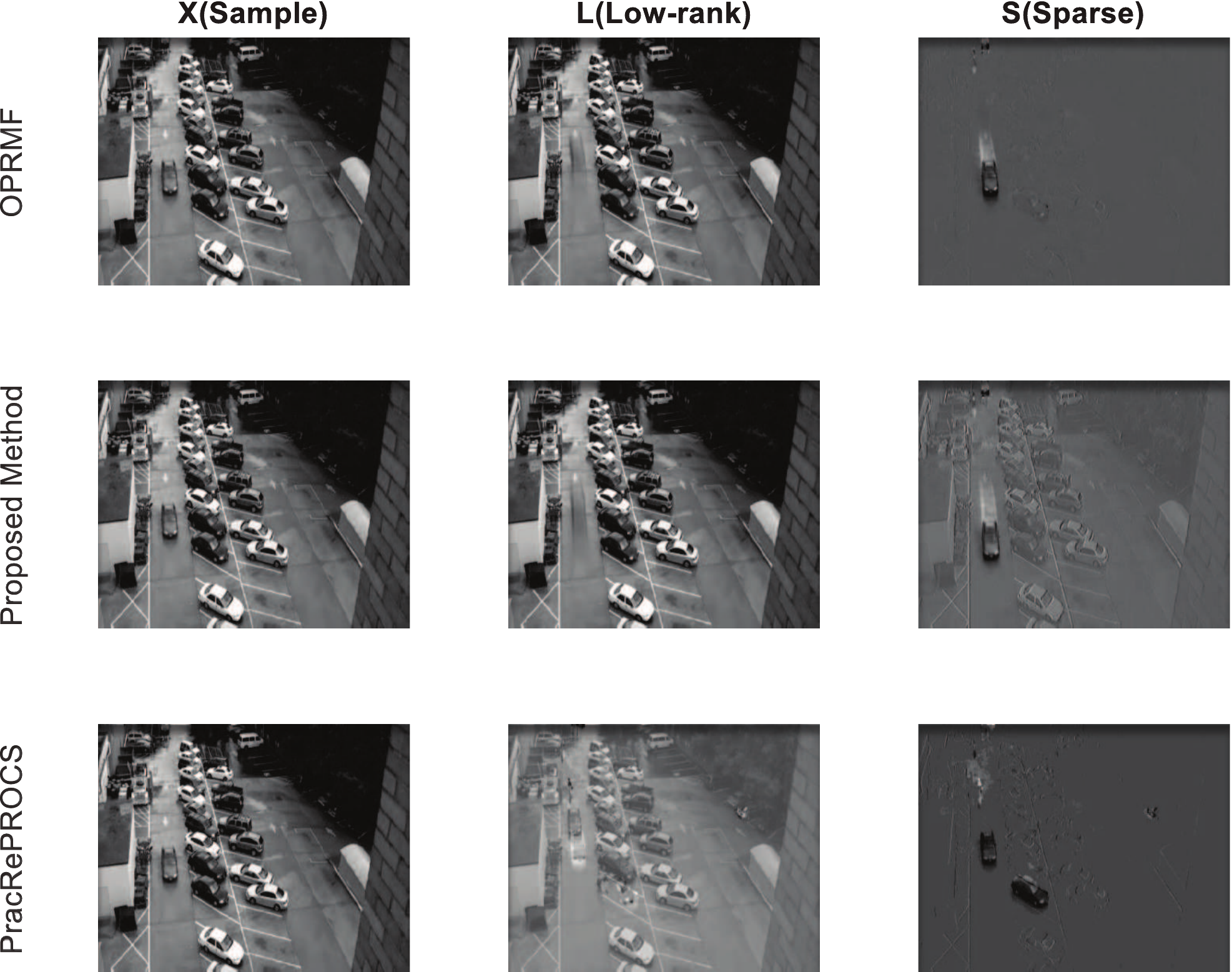}
	\caption{Performance comparisons with OPRMF \cite{wang2012probabilistic} and ReProCS \cite{reprocs} over dynamic background. The three algorithms are arranged, from top to bottom, in decreasing order of runtime. 
	}
	\label{dynamic}
\end{figure}

For the purposes of comparison, we consider two state-of-the-art foreground-background separation algorithms, namely OPRMF \cite{wang2012probabilistic} and ReProCS \cite{reprocs}. To keep the comparisons fair, we do not pre-train any of the algorithms and run them directly in an online manner while ignoring the first few frames. For all videos, we use $\mu_{\L}=0.005$, $\mu_{\S}=2$, $\lambda_{\L}=100$,  $\lambda_{\S}=.034$, and $\alpha_{\L}=\alpha_{\S}=0.2 $. All the experiments were  carried out on an Intel core i7-4770 computer running at 3.40 GHz with 16GB DDR3 3.40 GHz RAM. A 64-bit version of MATLAB 14 was used and execution times were measured using the default system clock.

Fig. \ref{static} considers a simple scenario with static background which is relatively easy to extract. However lack of training data results in ReProCS mis-classifying part of the foreground as background. Nevertheless, in such a setting, the IP-OGD algorithm does not outperform either of the two algorithms. 

Next, we consider the more complicated case of moving background arising in the video of a parking lot taken from \cite{virat}. Specifically, cars that are being parked are part of the foreground, while cars that are already parked constitute the background. As can be seen from Fig. \ref{dynamic}, ReProCS cannot deal with the dynamic background and continues to treat the parked car as foreground even after it has stopped. On the other hand, the OPRMF is relatively more accurate and only a small part of the parked car is included in the foreground. The proposed algorithm is slightly worse than OPRMF but has a significantly lower runtime \footnote{The readers are referred to the full video of the tests available at \url{https://bit.ly/2IkNbLz}.}.

The full gradient version of the algorithm is computationally expensive since it requires carrying out a singular value decomposition of the full matrix $\M_k$ at each $k$. The large-scale learning algorithms proposed in Sec. \ref{vrsec} can however be used to trade-off speed with accuracy. Towards this end, we reformulate the problem by defining $f_k^j$ as 
\begin{align}\label{here}
\!\!\!f_k^j = \norm{\Ob_j\odot(\M_k-\L-\S)}_F^2 + \mu_{\L}\norm{\L}_F^2 + \mu_{\S}\norm{\S}_F^2, 
\end{align}
for all $j = 1, \ldots, N$. The selection matrix $\Ob_j$ in \eqref{here} is defined according to the percentage of $\M_k$ sampled. For instance, with $12.5\%$ sampling, there are a total of $N = 8$ components, each depending on a block of size $r_s \times L$ where $r_s = r/8$, and
\begin{align}
\Ob_j = \begin{bmatrix} \mathbf{0}_{r_s(j-1) \times L} \\
\mathbf{1}_{r_s\times L} \\
\mathbf{0}_{r-jr_s\times L} 
\end{bmatrix}.
\end{align}
The results are presented for OP-SVRG algorithm with circular and randomly selected component index $i$ at each $k$. Observe that while with 25\% sampling the proposed algorithm runs four times faster, the loss in performance as compared to the full gradient case is very small.

\begin{figure}
	\centering
	\includegraphics[scale=0.8]{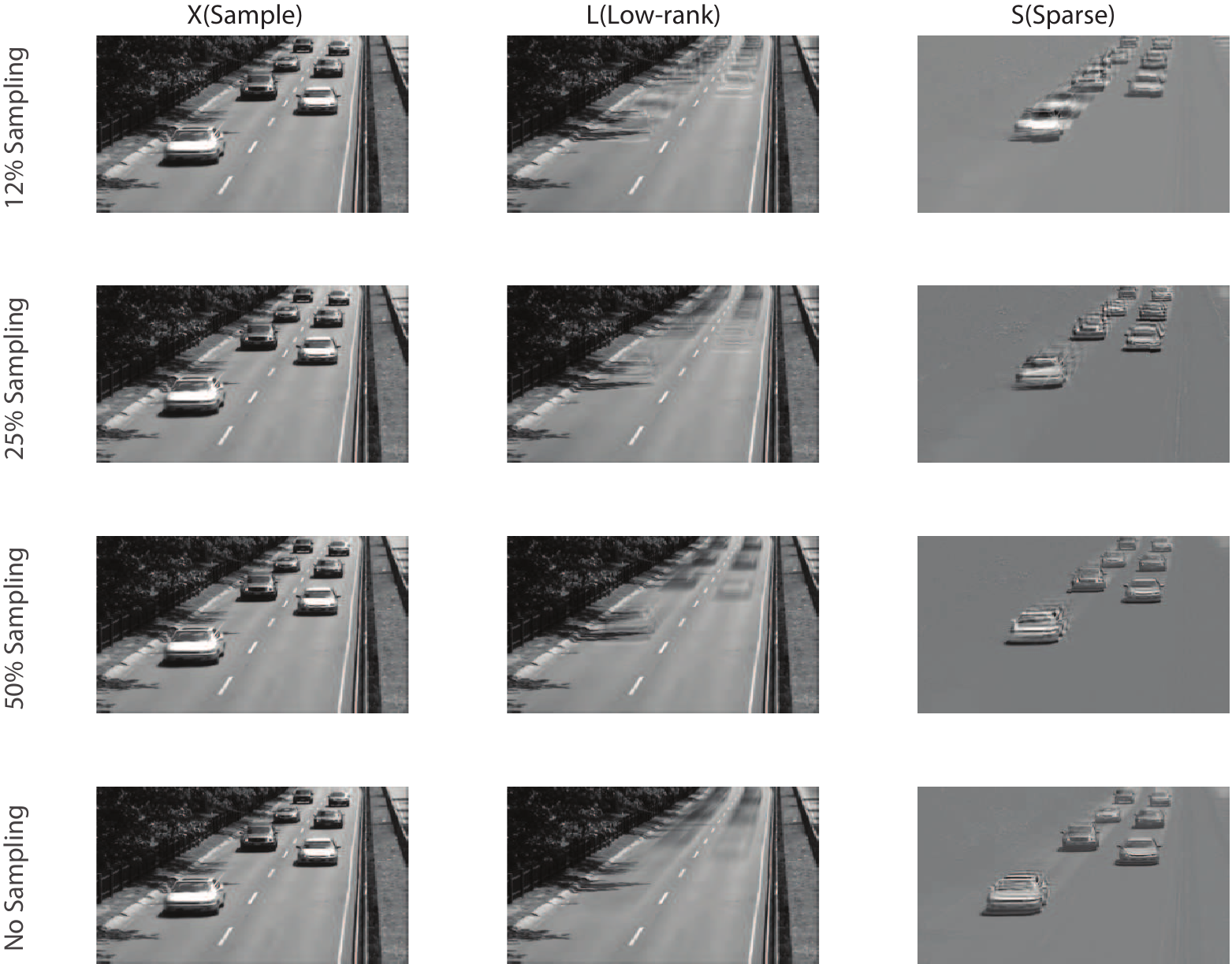}
	\caption{Separation of Low rank and sparse components using Proximal-SVRG under circular samplings}\label{circle_sampl}
\end{figure}

\begin{figure}
	\centering
	\includegraphics[scale=0.6]{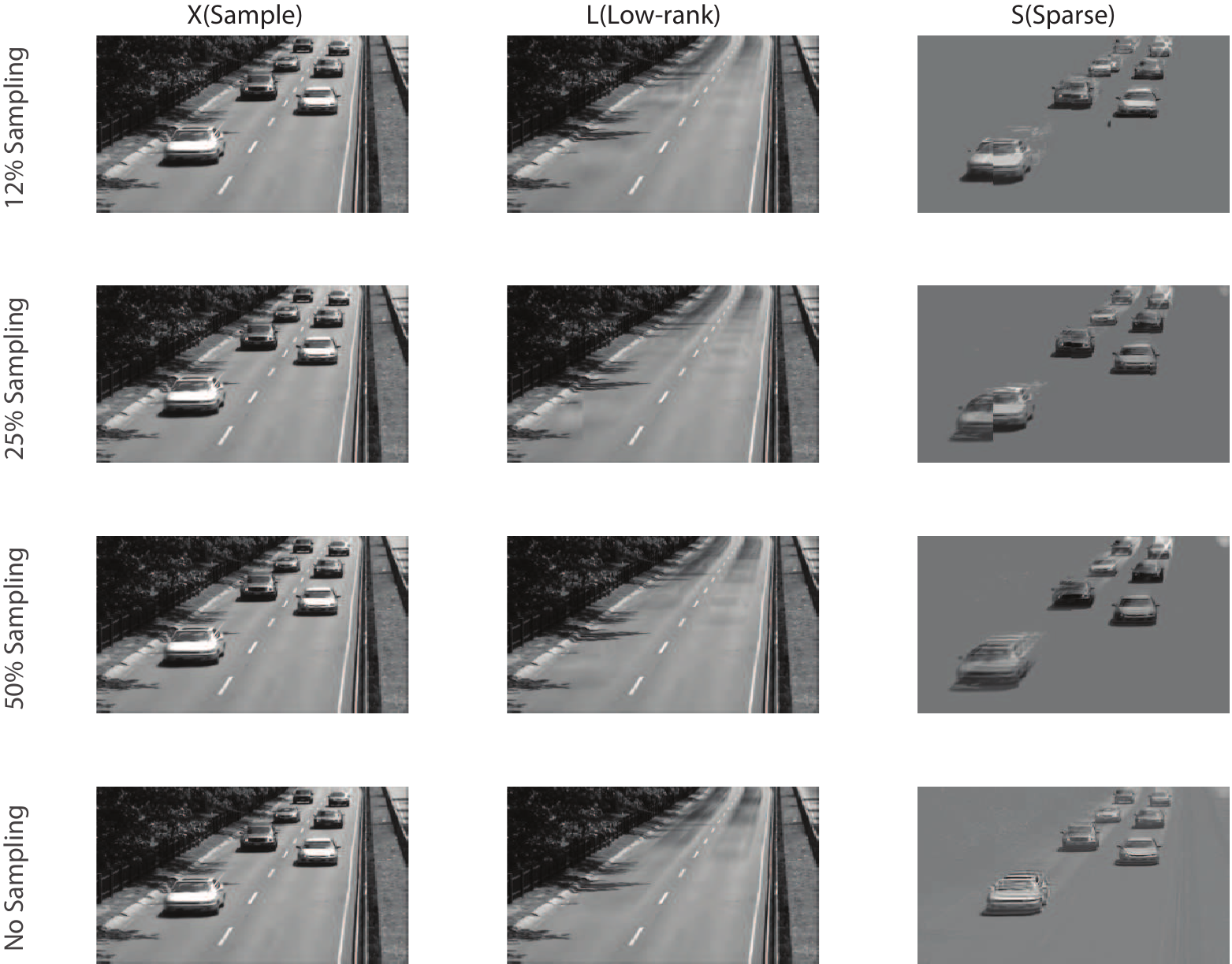}
	\caption{Separation of Low rank and sparse components using Proximal-SVRG under random samplings}\label{random_sampl}
\end{figure}

\section{Conclusion}  		 
This work focuses on the online convex optimization framework comprising of a sequential game between a learner and an adversary. A particular case of time varying optimization is considered when the given loss function is of the form of sum of a differentiable and non differentiable functions. This formulation is motivated considering the applications of \emph{robust subspace tracking} and \emph{bot formation trajectories}. A proximal inexact online gradient descent (OGD) method is proposed to solve the problem that tolerates changes in the function optimum as well as errors in the gradient. It is established that the proposed online algorithm exhibits sublinear dynamic regret provided that path length and cumulative gradient error are sublinear. 

For tackling large-scale problems, we proposed two different variants of the accelerated proximal OGD that do not require calculation of the full gradient function at every time instant. These algorithms can be viewed as the dynamic variants of the stochastic variance reduced gradient descent and incremental hybrid stochastic gradient methods. The performance of proposed algorithms is validated on a robot formation control problem and a robust subspace tracking problem applied to the dynamic foreground-background separation problem. 
\appendices
  		  
\section{Proofs of Lemmas \ref{lem1}-\ref{lemsg}}\label{proof:lemma1}
\subsection{Preliminaries}
Before discussing the proofs, we state some preliminary results that will be used in all the proofs. The Peter-Paul inequality for positive numbers $x$, $y$, and $c$ takes the form $(x+y)^2 \leq (1+c)x^2 + (1+1/c)y^2$. Since $f_k$ is strongly convex, the optimum $\xks$ is unique and satisfies the following two properties:
\begin{align}
&\xks = \pk{\xks - \alpha \nabla f_k(\xks)}, \label{xstar}\\
&\ip{\nabla f_{k}(\xks) + \boldsymbol{\nu},\x-\xks} \geq 0, \label{xopt}
\end{align} 
for all $\x \in \cX$, where $\boldsymbol{\nu} \in\partial g_{k}(\xks)$. 

Consider the proximal operator defined in \eqref{proximal} and let $\y := \pk{\w}$ for some $\w \in \Rn^n$. Then we have that
\begin{align}
\norm{\pk{\u} - \pk{\v}} &\leq \norm{\u-\v}  & \u,\v \in \cX\label{p1}\\ 
 \w-\y &\in \alpha \hspace{1pt} \partial  g_k(\y). \label{p2}
\end{align}

\begin{IEEEproof}[Proof of Lemma \ref{lem1}]
We begin with considering the difference between the next iterate $\x_{k+1}$ and the current optimal $\xks$ and using the properties \eqref{xstar} and \eqref{p1} as follows:
\begin{align}
&\norm{\x_{k+1}-\xks}^2 \nonumber\\
&=\norm{\pk{\x_k-\alpha\nt f_k(\x_k)} - \pk{\xks-\alpha\nabla f_k(\xks)}}^2  \nonumber\\
&\leq\norm{\x_k-\alpha\nt f_k(\x_k) - \xks + \alpha\nabla f_k(\xks)}^2 \label{eq:lemma_4_1}\\
&= \norm{(\x_k-\alpha\nabla f_k(\x_k)) - (\xks-\alpha\nabla f_k(\xks)) - \alpha\e_k}^2\label{eq:lemma_4_2}  \\
&= \norm{(\x_k-\alpha\nabla f_k(\x_k)) - (\xks-\alpha\nabla f_k(\xks))}^2 + \alpha^2\norm{\e_k}^2 \nonumber\\
&-2\alpha{\e_k}^{T}{[(\x_k-\alpha\nabla f_k(\x_k)) - (\xks-\alpha\nabla f_k(\xks))]},\label{eq:lemma_4_3} 
\end{align}
where the last equality in \eqref{eq:lemma_4_3} is obtained by expanding the square. The first term on the right hand side of \eqref{eq:lemma_4_3} can be expanded as 
\begin{align}	  		  
&\norm{(\x_k-\alpha\nabla f_k(\x_k)) - (\xks-\alpha\nabla f_k(\xks))}^2 \nonumber\\
&= \norm{\x_k-\xks}^2 + \alpha^2\norm{\nabla f_k(\x_k) - \nabla f_k(\xks)}^2 \nonumber\\
& - 2\alpha \ip{\nabla f_k(\x_k) - \nabla f_k(\xks), \x_k-\xks} \nonumber\\
&\leq  \Vert {\x}_{k} - \xks \Vert^{2} + \alpha^{2} L^{2} \Vert {\x}_{k} - \xks \Vert^{2} - 2 \alpha \mu \Vert {\x}_{k} - \xks \Vert^{2}  \label{eq:lemma_4_35} \\
 &=   \rho^2 \Vert {\x}_{k} - \xks \Vert^{2} \label{eq:lemma_4_4} ,
	  	\end{align}
	  				where, $\rho^2:=(1 + \alpha^{2}L^{2} - 2 \alpha\mu )$ and $\alpha$ is chosen such that $\rho<1$. The second inequality in \eqref{eq:lemma_4_35} holds since the function $f_k$ is strongly convex with parameter $\mu$ and $\nabla f_k$ is Lipschitz continuous with parameter $L$. Taking the positive square root on both sides of \eqref{eq:lemma_4_4}, we also have that $\norm{(\x_k-\alpha\nabla f_k(\x_k)) - (\xks-\alpha\nabla f_k(\xks))} \leq \rho \norm{\x_k-\xks}$. Substituting \eqref{eq:lemma_4_4} into \eqref{eq:lemma_4_3}, we obtain 	  				
	  				\begin{align}
	  			 \norm{{{\x}}_{k+1}\!-\!\xks}^{2} \!\! \leq &  \rho^{2} \norm{{\x}_{k} \!-\! \xks}^{2}  \!\!+\! \alpha^2\norm{\e_{k}}^{2}\!\!\! + \!\!2\alpha\rho \norm{\e_{k}} \norm{\bar{\x}_{k}\!-\!\xks} \nonumber\\
	  	\leq & (\rho \norm{{\x}_{k}-\xks} + \alpha\norm{\e_k})^{2}  \label{eq:lemma_4_5},
	  		  				\end{align}
	  				with $\delta_k:=\norm{\alpha\e_k}$. The required result follows from taking the positive square root on both sides of \eqref{eq:lemma_4_5}. 
	  		\end{IEEEproof}
	
				
 \begin{IEEEproof}[Proof of Lemma \ref{lem2}] 
We begin with considering cumulative optimality gap and applying the triangle inequality as follows. 
\begin{align}
\sum_{k=2}^K\norm{\x_k- \xks} = & \sum_{k=2}^K\norm{\x_k- \x_{k-1}^\star + \x_{k-1}^\star - \xks} \nonumber\\
	  		 	 \leq & \sum_{k=2}^K\norm{\x_k- \x_{k-1}^\star} +  \sum\limits_{k=2}^K\norm{\xks - \x_{k-1}^\star} \nonumber \\
= & \sum\limits_{k=1}^{K-1}\norm{\x_{k+1}- \xks} + W_K\label{eq:corr1}	,
\end{align}
where the equality in \eqref{eq:corr1} follows from the definition of the path length. Next, utilizing the result of Lemma 4, we obtain
\begin{align}
\sum_{k=2}^{K}\norm{\x_{k}- \xks } &\leq  \sum_{k=1}^{K-1}(\rho \norm{\x_{k}- \xks} + \alpha\norm{\e_k}) + W_{K} \label{eq:complete_K} \\
&\leq \sum_{k=1}^K(\rho \norm{\x_{k}- \xks} + \alpha\norm{\e_k}) + W_{K} \label{eq:complete_K2},
\end{align}
where the inequality in \eqref{eq:complete_K2} holds due to the inclusion of an additional non-negative terms on the right. Adding $\norm{\x_1-\x_1^\star}$ on both sides of 
\eqref{eq:complete_K2}, we obtain the required result as follows.  
	  		    \begin{align}	\label{eq:complete_K_2}
	  		    &\sum_{k=1}^{K}\norm{\x_{k}- \xks } \nonumber\\
						&\leq   \norm{\x_{1} - \x_{1}^\star} + \sum_{k=1}^{K}\rho \norm{\x_{k}- \xks} + \alpha\sum\limits_{k=1}^K \norm{\e_k}  + W_K \\
						&= \frac{1}{1-\rho}\norm{\x_{1} - \x_{1}^\star}  + \frac{\alpha}{1-\rho}\sum\limits_{k=1}^K \norm{\e_k}  + \frac{1}{1-\rho}W_K,
	  		    \end{align}
where the last equality follows simply from taking the cumulative optimality gap to the right. 
 \end{IEEEproof}		  		

\begin{IEEEproof}[Proof of Lemma \ref{lemsg}]	
Recalling the proximal update equation $\x_{k+1} = \pk{\x_k-\alpha\nt f_k(\x_k)}$ and using the property \eqref{p2}, we obtain
\begin{align}
\x_{k}-\alpha\nt f_k(\x_k) -\x_{k+1} &\in \alpha \partial g_k(\x_{k+1}) ,
\end{align}
since we have $\x_{k+1}\in\cX$ which implies that $\boldsymbol{\mathcal{N}}_{\cX}(\x_{k+1}) = \boldsymbol{0}$ where $\boldsymbol{\mathcal{N}}_{\cX}(\x_{k+1})$ is the norm cone at $\x_{k+1}$ w.r.t. convex set $\cX$.
Equivalently, defining $G_k(\x_k):=\nt f_k(\x_k) + \u_k$ where $\u_k \in \partial g_k(\x_{k+1})$, the update equation can also be written as
\begin{align}\label{update}
\alpha G_k(\x_{k}) &= \x_{k} - \x_{k+1} \\
\Rightarrow \alpha G_k(\x_{k}) = & (\x_{k} - \xks) - (\x_{k+1} - \xks).
\end{align}
Taking norm on both sides and applying the triangle inequality on the right, we obtain\begin{align}
\alpha \norm{G_k(\x_{k})} &\leq  \norm{\x_k - \xks} + \norm{\x_{k+1} - \xks} \\
&\leq \norm{\x_{k} - \xks} + \rho \norm{\x_{k} - \xks}  +\alpha\norm{\e_k} \label{b0} \\
\Rightarrow \norm{G_k(\x_{k})} &\leq \frac{1+\rho}{\alpha} \norm{\x_{k} - \xks}    +\norm{\e_k}  \label{b1},
\end{align}
where we have used the result from Lemma \ref{lem1} in \eqref{b0}. In order to bound $\norm{h_k(\x_k)}$ however, we need to express $\partial h_k(\x_{k})$	in terms of $G_k(\x_{k})$ as follows:
\begin{align}
\partial  h_k(\x_{k}) &=  \nabla f_{k}(\x_{k}) + \partial g_{k}(\x_{k}) \\
&\hspace{-1cm}=  \nabla f_{k}(\x_{k}) + \partial g_{k}(\x_{k+1}) + \partial g_{k}(\x_{k}) - \partial g_{k}(\x_{k+1}) \\
&=  G_{k}(\x_{k}) - \e_k  + \partial g_{k}(\x_{k}) - \partial g_{k}(\x_{k+1}) .\label{hkgk}
\end{align}	 
Next, taking norm on both sides of \eqref{hkgk} and using the triangle inequality, we obtain
	  			\begin{align}
	  			\norm{\partial  h_k(\x_{k})} &\leq  \norm{G_{k}(\x_{k})} + \norm{\e_k}  + \norm{\partial g_{k}(\x_{k})} + \norm{\partial g_{k}(\x_{k+1})} \\
					&\leq \frac{1+\rho}{\alpha} \norm{\x_{k} - \xks}  + 2\norm{\e_k}  + 2L_g \label{sbound},
	  			\end{align}
where \eqref{sbound} follows from the use of \eqref{hkgk} and the bounded subgradient property that follows from the Lipschitz continuity  property of $g_k$ in \eqref{Lip_g1}.

Next, we upper bound the term $\EE\norm{\x_{k} - \xks}$ as follows:
\begin{align}
\EE\norm{\x_{k} - \xks} &=  \EE\norm{\x_{k} - \x_{k-1}^\star + \x_{k-1}^\star - \xks} \\
	&\leq  \EE\norm{\x_{k} - \x_{k-1}^\star} + \norm{\xks - \x_{k-1}^\star}\\
	&\leq  \rho \EE\norm{\x_{k-1} - \x_{k-1}^\star} + \alpha \EE\norm{\e_{k-1}} + \sigma \label{xkbnd0}\\
	&\leq \rho \EE\norm{\x_{k-1} - \x_{k-1}^\star} + \alpha \gamma + \sigma	\label{xkbnd},			
\end{align}
where the inequality in \eqref{xkbnd} follows from the result of Lemma \ref{lem1} and the bounded variations property \eqref{eq:BV}, while that in \eqref{xkbnd} follows from \eqref{eq:BV2}. The bound in \eqref{xkbnd} can be recursively applied to yield
\begin{align}
\EE\norm{\x_k - \xks}
&\leq  \rho \EE\norm{\x_{k-1} - \x_{k-1}^\star} + \alpha \gamma + \sigma  \nonumber\\
&\!\!\!\!\!\!\!\leq  \rho^{k-1} \norm{\x_{1} - \x_{1}^\star} + (\alpha \gamma + \sigma)\bigg( \frac{1 - \rho^{k-1}}{1-\rho} \bigg)\label{xkbndf0} \\
\leq &  \norm{\x_{1} - \x_{1}^\star} + \frac{\alpha \gamma + \sigma}{1-\rho} \label{xkbndf}  .
\end{align}
Finally, taking expectation in \eqref{sbound} and using the bound in \eqref{xkbndf}, we obtain
\begin{align}
\!\!\!\!\EE\norm{\partial  h_k\!(\!\x_k\!)\!} \!\leq  \! \frac{1\!\!+\!\!\rho}{\alpha} \bigg(\!\!\!\norm{\x_1-\x_1^\star} +  \frac{\alpha \gamma + \sigma}{1-\rho} \bigg)  +   2\gamma + 2L_g,
\end{align}
which is the required result. \end{IEEEproof}

\section{Proof of Lemma \ref{lemsvrg}} \label{proof_pro2}	
Observe that Assumptions \textbf{(A1')-(A2')} also imply \textbf{(A1)-(A2)} for the full function $f_k$. Consequently, the results in Lemma 1 and 2 hold directly. Moreover, Lemma 3 follows with $D = L_g + M$. Let $\iota_k$ be the random index that is selected at time $k$ and let $\kappa = \lfloor k/K_0\rfloor K_0 + 1$ be the iteration when the entire memory was selected previously. It can be seen that $0\leq k-\kappa \leq K_0-1$. For the sake of brevity, for any mapping $\psib:\cX \rightarrow \Rn^n$, let $\Theta(\psib):=\sup_{\x\in\cX} \norm{\psib(\x)}$.

When analyzing the error at iteration $k$, we will drop the subscript and denote the component index by $\iota$. Similarly, we denote the expectation with respect to $\iota$ as $\EE_\iota[\cdot]$. As $\iota_k$ are independent identically, distributed, $\EE_{\iota}[\cdot]$ is equivalent to the conditional expectation $\EE[\cdot \mid \iota_1, \ldots, \iota_{k-1}]$. Further since $\iota$ is uniformly distributed between 1  and $N$, we have that
\begin{align}
\EE_{\iota}\left[\nabla f_k^\iota(\x_k)\right] &= \frac{1}{N}\sum_{j=1}^N \nabla f_k^j(\x_k) = \nabla f_k(\x_k),\\
\EE_{\iota}[\bt^\iota] &= \frac{1}{N}\sum_{j=1}^N \bt^j .
\end{align}
These identities imply that $\e_k$ is zero-mean and the mean-squared error can be written as
\begin{align}
&\EE_{\iota}\norm{\e_k}^2 = \EE_{\iota}\norm{\nabla f_k^{\iota}(\x_k)-\bt^{\iota} - \EE_{\iota}[\nabla f_k^{\iota}(\x_k)-\bt^{\iota}]}^2 \\
&=\EE_{\iota}\norm{\nabla f_k^{\iota}(\x_k)-\bt^{\iota}}^2 -\norm{\EE_{\iota}[\nabla f_k^{\iota}(\x_k)-\bt^{\iota}]}^2 \label{compacterror},
\end{align}
where $\beta^\iota = \nabla f_\kappa^\iota(\x_\kappa)$. Adding and subtracting $f_k^\iota(\x_\kappa)$ within the first term of \eqref{compacterror} and applying the Peter-Paul inequality with some $c > 0$, we obtain
\begin{align}
&\norm{\nabla f_k^{\iota}(\x_k)-\nabla f_{\kappa}^{\iota}(\x_{\kappa})}^2 
 \leq  (1+ c)\norm{\nabla f_k^{\iota}(\x_k)-\nabla f_k^{\iota}(\x_{\kappa})}^2\nonumber \\ &\hspace{2cm}+(1+\tfrac{1}{c})\norm{\nabla f_k^{\iota}(\x_{\kappa})-\nabla f_{\kappa}^{\iota}(\x_{\kappa})}^2 \nonumber \\
& \leq \!\! (1\!+\!c)L^2\norm{\x_k\!-\!\x_{\kappa}}^2\!\!\!+\!(1\!\!+\!\!\tfrac{1}{c})\norm{\nabla f_k^{\iota}(\x_{\kappa})\!-\!\nabla f_{\kappa}^{\iota}(\x_{\kappa})}^2  \label{firsterrorbound} 
\end{align}
where we have used the fact that $\nabla f_k^\iota$ is Lipschitz continuous. The constant $c$ is left unspecified at this point and will be fixed later. Similarly, in order to bound the second term in \eqref{compacterror}, we again consider the Peter-Paul inequality
\begin{multline}
\norm{\EE_{\iota}[\nabla f_k^{\iota}(\x_k)-\nabla f_{k}^{\iota}(\x_{\kappa})]}^2 \\
\leq 	(1+a)\norm{\EE_{\iota}[\nabla f_k^{\iota}(\x_k)-\nabla f_{\kappa}^{\iota}(\x_{\kappa})]}^2  \\
+(1+\tfrac{1}{a})\norm{\EE_{\iota}[\nabla f_{\kappa}^{\iota}(\x_{\kappa})-\nabla f_k^{\iota}(\x_{\kappa})]}^2.  
\end{multline}
for some $a>0$. 
Dividing both sides by $(1+a)$ and rearranging, we obtain,
\begin{align}
&-\norm{\EE_{\iota}[\nabla f_k^{\iota}(\x_k)-\nabla f_{\kappa}^{\iota}(\x_{\kappa})]}^2 \nonumber\\
&\hspace{1cm} \leq  -\tfrac{1}{1+a}\norm{\EE_{\iota}[\nabla f_k^{\iota}(\x_k)-\nabla f_{k}^{\iota}(\x_{\kappa})]}^2 \nonumber \\  
&\hspace{2cm} + \tfrac{1}{a}\norm{\EE_{\iota}[\nabla f_{\kappa}^{\iota}(\x_{\kappa})-\nabla f_k^{\iota}(\x_{\kappa})]}^2 \nonumber \\
&\hspace{-2mm}\leq-\tfrac{\mu^2}{1+a}\norm{\x_k-\x_{\kappa}}^2  + \tfrac{1}{a}\norm{\EE_{\iota}[\nabla f_{\kappa}^{\iota}(\x_{\kappa})-\nabla f_k^{\iota}(\x_{\kappa})]}^2\label{seconderrorbound},
\end{align}
where \eqref{seconderrorbound} follows from the strong convexity of $f_k^\iota$. The last terms in \eqref{firsterrorbound} and \eqref{seconderrorbound} can both be bounded by observing that 
\begin{align}
\norm{\EE_{\iota}[\nabla f_{\kappa}^{\iota}(\x_{\kappa})-\nabla f_k^{\iota}(\x_{\kappa})]}^2 &\leq \EE_\iota\norm{\nabla f_{\kappa}^{\iota}(\x_{\kappa})-\nabla f_k^{\iota}(\x_{\kappa})}^2 \nonumber\\
&\hspace{-2cm}\leq \frac{1}{N}\sum_{j=1}^N \norm{\nabla f_{\kappa}^{\iota}(\x_{\kappa})-\nabla f_k^{\iota}(\x_{\kappa})}^2\label{iidiota}\\
&\hspace{-2cm}\leq \frac{1}{N}\sum_{j=1}^N\Theta^2(\nabla f_\kappa^j-\nabla f_k^j)\label{theta},
\end{align}
where \eqref{iidiota} follows since the index $\iota$ is chosen independently and with equal probability. Thus, the final bound on the mean-squared error becomes
\begin{align}
\EE_{\iota}\norm{\e_k}^2 &\leq  \Gamma^2\norm{\x_k-\x_{\kappa}}^2 + \zeta^2\sum_{j=1}^N\Theta^2(\nabla f_{\kappa}^j-\nabla f_k^j) \label{lipsaga},
\end{align}
where $\Gamma^2 :=  L^2(1+c)-\tfrac{\mu^2}{(1+a)}$ and $\zeta^2:=\frac{1}{N}(1+\tfrac{1}{c}+\tfrac{1}{a})$. 
Equivalently, from the Cauchy-Schwarz inequality, we have that  
\begin{align}
\EE_\iota\!\!\norm{\e_k} \!\!\leq\!\! \sqrt{\!\EE_{\iota}\!\norm{\e_k}^2}\!\!\leq\!\! \Gamma\!\norm{\x_k\!\!-\!\!\x_{\kappa}} \!\!+\!\!\zeta\sum_{j=1}^N\Theta(\nabla f_k^j - \nabla f_{\kappa}^j)\label{svrgbnd1}.
\end{align}
For every iteration $k$, the previous time instant at which entire memory was updated is a unique index $\kappa = \lfloor (k-1)/K_0\rfloor K_0+1$ that depends on $k$ but is constant for all $\kappa \leq k < \kappa+K_0$. Therefore we have that 
\begin{align}
 \sum_{k=1}^K\norm{\x_{k}-\x_{\kappa}}&\leq \sum_{k=1}^K(\norm{\x_{k}-\x_{k}^*}+ \norm{\x_{\kappa}-\x_{\kappa}^*}+ \norm{\x_{k}^*-\x_{\kappa}^*})  \nonumber\\
  &  \leq \sum_{k=1}^K\norm{\x_{k}-\x_{k}^\star} + \sum_{k=1}^{K}\norm{\x_k^\star-\x_{\kappa}^\star}\nonumber\\
	&\hspace{-1cm}+ (K_0-1)\sum_{j=0}^{\lfloor (K-1)/K_0 \rfloor}\norm{\x_{jK_0+1}-\x_{jK_0+1}^\star}. \label{randompath} 
\end{align}
where the right hand side of \eqref{randompath} included additional non-negative terms. The first term in \eqref{randompath} is bounded from Lemma \ref{lem2}. The second term in \eqref{randompath} can be bounded using the triangle inequality and inclusion of additional terms as follows
\begin{align}
\sum_{k=1}^K\norm{\x_k^\star-\x_\kappa^\star} &\leq \sum_{k=1}^K\sum_{r=0}^{\min\{k-2,K_0-1\}}\norm{\x_{k-r}^\star-\x_{k-r-1}^\star} \nonumber\\
&\leq K_0 \sum_{k=1}^K \norm{\x_k^\star-\x_{k-1}^\star} = K_0W_K.
\end{align}
Note that the second term in \eqref{svrgbnd1} can also be bounded along similar lines
\begin{align}
\sum_{k=1}^K \Theta(\nabla f_k^j - \nabla f_{\kappa}^j) &\leq K_0\sum_{k=1}^K \Theta(\nabla f_k^j - \nabla f_{k-1}^j) \\
&\leq K_0\sqrt{KV_K^j}.
\end{align}

In order to obtain a bound on the third term in \eqref{randompath}, consider 
\begin{align}
&\EE\norm{\x_{jK_0+1}\!\!-\!\!\x_{jK_0+1}^\star}\!\!\leq\!\! \EE\norm{\x_{jK_0+1}\!\!-\!\!\x_{jK_0}^\star} \!\!+\!\! \EE\norm{\x_{jK_0+1}^\star\!\!-\!\!\x_{jK_0}^\star} \nonumber\\
&\leq \rho\EE\norm{\x_{jK_0}-\x_{jK_0}^\star} + \alpha \EE\norm{\e_{jK_0}} + \norm{\x_{jK_0+1}^\star-\x_{jK_0}^\star}, \label{jcontraction}
\end{align}
for any $j \geq 0$, where \eqref{jcontraction} follows from Lemma \ref{lem1}. Recursively applying the inequality in \eqref{jcontraction}, we have that
\begin{align}
&\EE\norm{\x_{jK_0+1}-\x_{jK_0+1}^\star} \leq \rho^{K_0}\EE\norm{\x_{(j-1)K_0+1}-\x_{(j-1)K_0+1}^\star} \nonumber\\
&\hspace{1.2cm}+ \alpha \sum_{r = 1}^{K_0}\rho^{K_0-r}\EE\norm{\e_{(j-1)K_0+r}} \nonumber\\
&\hspace{1.2cm}+  \sum_{r = 1}^{K_0}\rho^{K_0-r}\norm{\x_{(j-1)K_0+r}^\star-\x_{(j-1)K_0+r-1}^\star}.
\end{align}
Taking summation over all $0\leq j \leq \lfloor(K-1)/K_0\rfloor$ and rearranging, we obtain
\begin{align}
(1-\rho^{K_0})\!\!&\sum_{j=0}^{\lfloor \frac{(K-1)}{K_0} \rfloor}\EE\norm{\x_{jK_0+1}-\x_{jK_0+1}^\star}\nonumber\\
&\hspace{-1cm}\leq \norm{\x_1-\x_1^\star} + \sum_{j=0}^{\lfloor \frac{(K-1)}{K_0} \rfloor}\sum_{r=1}^{K_0}\rho^{K_0-r}\EE\norm{\e_{(j-1)K_0+r}}\nonumber\\
&\hspace{-1cm}+ \sum_{j=0}^{\lfloor \frac{(K-1)}{K_0} \rfloor}\sum_{r=1}^{K_0}\rho^{K_0-r} \norm{\x_{(j-1)K_0+r}^\star-\x_{(j-1)K_0+r-1}^\star} \nonumber\\
&\hspace{-1.3cm}\leq \norm{\x_1-\x_1^\star} + \!\!\sum_{k=1}^{K}\rho^{K_0-\text{mod}(k,K_0)} \left[\alpha\EE\norm{\e_k} + \norm{\x_k^\star-\x_{k-1}^\star}\right] \nonumber\\
&\hspace{-1.3cm}\leq \norm{\x_1-\x_1^\star} + \alpha E_K + W_K,
\end{align}
where the last inequality simply uses the fact that $\rho^{K_0-\text{mod}(k,K_0)} < 1$ for all $k\geq 1$. Taking expectation in \eqref{randompath} and using the bounds for the three terms, we obtain
\begin{align}
&\sum_{k=1}^K\EE\norm{\x_{k}-\x_{\kappa}}\label{xkkappa}\\
&\leq (\frac{K_0}{1-\rho^{K_0}} + \frac{1}{1-\rho})(\norm{\x_1-\x_1^\star}+\alpha E_K + W_K) + K_0W_K.  \nonumber
\end{align}
Next, taking full expectation in \eqref{svrgbnd1}, using \eqref{xkkappa}, and rearranging, we obtain
\begin{align}
E_K \leq &\eta (\norm{\x_1-\x_1^\star}+\alpha E_K + W_K) + \Gamma K_0W_K \nonumber\\
& +\zeta K_0\sum_{j=1}^N\sqrt{KV_K^j} \\
&\hspace{-1.3cm} \leq \frac{\norm{\x_1-\x_1^\star} + (\Gamma K_0+\eta)W_K + \zeta K_0\sum_{j=1}^N\sqrt{KV_K^j}}{1-\eta\alpha}\\
&\hspace{-1.3cm}\leq \mathcal{O}(1+ W_K + \sum_{j=1}^N\sqrt{KV_K^j}),
\end{align}
where $\eta := \Gamma (\frac{K_0}{1-\rho^{K_0}} + \frac{1}{1-\rho})$ and it is required that $\alpha\eta < 1$. Equivalently, it is required that 
\begin{align}
\alpha\sqrt{L^2(1+c)-\tfrac{\mu^2}{1+a}}\left(\frac{K_0}{1-\rho^{K_0}} + \frac{1}{1-\rho}\right) < 1.
\end{align}
It can be verified that such a choice of parameters $\alpha < 2\mu/L^2$,  $a > 0$, and $c > 0$ is possible for any given $K_0$ when $\mu/L > 0.89$.

\section{Proof of Proposition \ref{lemiss}}\label{proof_prop_3}
We begin by observing that $\e_k$ in \eqref{error_incresing} may be written as
\begin{align}
\e_k &= \frac{1}{\abs{\I_k}}\sum\limits_{i \in \I_k}\nabla f_k^i(\x_k)- \frac{1}{N}\sum\limits_{i=1}^{N}\nabla f_k^{i}(\x_k) \nonumber\\
&=\frac{N-\abs{\I_k}}{N\abs{\I_k}}\sum\limits_{i \in \I_k}\nabla f_k^i(\x_k)- \frac{1}{N}\sum\limits_{i\notin \I_k}\nabla f_k^{i}(\x_k).
\end{align}
Therefore we have the following bound
			\begin{align}
			\!\!\!\norm{\e_k} &\leq  \frac{N\!-\!\abs{\I_k}}{N\abs{\I_k}}\!\!\sum\limits_{i \in \I_k}\!\!\norm{\nabla f_k^i(\x_k)}+ \frac{1}{N}\sum\limits_{i\notin \I_k}\norm{\nabla f_k^{i}(\x_k)} \\
			&\leq  \frac{N-\abs{\I_k}}{N\abs{\I_k}}\sum\limits_{i \in \I_k}M+ \frac{1}{N}\sum\limits_{i\notin \I_k}M  \\
			&\leq  \frac{N-\abs{\I_k}}{N\abs{\I_k}}M\abs{\I_k}+ \frac{N-\abs{\I_k}}{N}M  \leq 2M\frac{N-\abs{\I_k}}{N} \nonumber
			\end{align}
where we have used the triangle inequality as well as the gradient boundedness property \textbf{(A3')}. It can be seen that in order to bound the cumulative error norm, it is required that the quantity $N-\abs{\I_k}$ be summable. For instance for the choice $\abs{\I}_k  = N-N\exp(-\beta k)$ leads to the bound
\begin{align}\nonumber
\sum_{k=1}^K\EE\norm{\e_k} \leq \sum_{k=1}^K 2M\exp(-\beta k) \leq \frac{2M}{1-\exp(-\beta)}-1
\end{align}
for any $\beta > 0$. Similarly, the choice $\abs{\I_k} = N-\tfrac{N}{k^2}$ also leads to bounded cumulative error.

  		\vspace{-0mm}
  		
\footnotesize
\bibliographystyle{IEEEtran} 
\bibliography{IEEEabrv,ref}
\newpage\onecolumn
\section*{Supplementary material for "Online Learning with Inexact Proximal Online Gradient Descent Algorithms"}
\section{Details of Bot Formulation Problem}\label{bot_form}
At the beginning of any time instant $k$, we have an initial pose $\p(k)\in\mathcal{F}(k)$ of the robots including the leader where $\p(k) \neq \s(k) $, and we need to obtain a pose $\p(k+1) \in \mathcal{F}(k)$ for the bots such that $\p(k+1) = \s(k) $. 
Consider the case when all the robots are present on a plane ($d=2$) with $m$ follower bots and $1$ leader bot. We can represent the shape $\s(k)=\text{vec}{(\H(k)^T)}$, where $\H(k)$ is defined as
\begin{align}\nonumber
\!\!\!\!\H(k)  = \lbrace \omega \B(k)\R + 1_{m+1}\t^{T}  : \omega \in \Rn_{+}, \R \in SO(2), \t \in \Rn^2 \rbrace.
\end{align} 
In this formulation, $\omega$ is the positive scaling factor, $\R$ is the rotation matrix of size $2\times2$ with angle of rotation $\theta$, $\t$ is the translation vector which moves the base shape $\B(k)$ of size $(m+1)\times2$ in $\Rn^2$. This base shape $\B(k)$ is centered around the origin of plane and may vary with $k$.  

Let $\s_i(k) = (\s_i^x(k), \s_i^y(k))$ and $\p_i(k+1) = (\p_i^x(k+1),\p_i^y(k+1))$ denote Cartesian coordinates in $\mathcal{F}(k)$. If we set first bot to be the leader, and choose $\t = \p_1(k+1)$, then we have the following non-linear equality constraints of the form 
\begin{align}
\p_i^x(k+1) - \p_1^x(k+1) &= \omega (\s_i^x(k) \cos\theta - \s_i^y(k) \sin \theta),\nonumber \\
\p_i^y(k+1) - \p_1^y(k+1) &= \omega (\s_i^x(k) \sin\theta + \s_i^y(k) \cos \theta) ,
\end{align}
for $i = 2, \cdots, m+1$. Without loss of generality, we can define the formation orientation $\theta$ and scale $\omega$ as :
\begin{align}
\theta = \arctan \frac{\p_2(k+1)^y - \p_1^y(k+1)}{\p_2^x(k+1) - \p_1^x(k+1)} \hspace{0.5cm} \text{and}  \hspace{0.5cm} \omega = \frac{\norm{\p_2(k+1) - \p_1(k+1)}_2}{\norm{\s_2(k)}_2},
\end{align}
where we have assumed that $\s_1(k) = O_{\mathcal{F}(k)}$. Using the above substitution, we can rewrite the constraints as 
\begin{align}
\norm{\s_2(k)}_2(\p_i^x(k+1) - \p_1^x(k+1)) &= (\s_i^x(k),-\s_i^y(k))^T (\p_2(k+1) - \p_1(k+1)), \nonumber \\
\norm{\s_2(k)}_2(\p_i^y(k+1) - \p_1^y(k+1)) &= (\s_i^y(k),-\s_i^x(k))^T (\p_2(k+1) - \p_1(k+1)),
\end{align}
for $i = 3,\cdots,m+1$. Now we are left with $2\times(m-1)$ linear functions of our state vector $\p(k+1) = [\p_1(k+1), \p_2(k+1), \cdots, \p_{m+1}(k+1)]^T$ which we need to estimate. These constraints sufficiently encode the desired shape of the formation. Moreover, these constraints can be compactly written as $\A(k) \p(k+1) = \textbf{0}$ where $\A(k)$ is some low rank matrix.

\end{document}